\documentclass[11pt]{article}

\usepackage{graphics}
\usepackage{amsfonts}
\usepackage{amssymb}
\usepackage{amsmath}
\usepackage{amsthm}
\usepackage{url}
\usepackage{color}
\usepackage{multirow}
\usepackage{enumerate}
\usepackage{hyperref}

\newtheorem{thm}{Theorem}
\newtheorem{lem}{Lemma}
\newtheorem{cor}{Corollary}
\newtheorem{prop}{Proposition}

\theoremstyle{definition}%an alternative is \theoremstyle{remark}
\newtheorem{rem}{Remark}

\theoremstyle{definition}
\newtheorem{deff}{Definition}

\theoremstyle{definition}
\newtheorem{example}{Example}

\newcommand{\R}{\mathbb{R}}

\newcommand{\cA}{\mathcal{A}}

\newcommand{\cS}{\mathcal{S}}
\DeclareMathOperator{\tr}{tr}
\DeclareMathOperator{\diag}{diag}
\DeclareMathOperator{\Diag}{Diag}
\DeclareMathOperator{\aut}{aut}
\DeclareMathOperator{\spen}{span}
\newcommand{\vect}{\mbox{vec}}
\newcommand{\Proof}{{\em Proof}.~}
\renewcommand{\binom}[2]{\genfrac(){0pt}{}{#1}{#2}}

\title{On bounding the bandwidth of graphs with symmetry\footnote{This version is published in INFORMS Journal on Computing 27 (2015), 75--88.}}
\author{ {E.R. van Dam}\thanks{Department of Econometrics and OR, Tilburg
University, The Netherlands. {\tt edwin.vandam@uvt.nl} }
 \and {R. Sotirov}\thanks{Department of Econometrics and OR, Tilburg
University, The Netherlands. {\tt r.sotirov@uvt.nl} }}

\date{}

\begin{document}
\maketitle

\begin{abstract}
We derive a new lower bound for the bandwidth of a graph that is based on a new
lower bound for the minimum cut problem. Our new semidefinite programming
relaxation of the minimum cut problem is obtained by strengthening the known
semidefinite programming relaxation for the quadratic assignment problem (or
for the graph partition problem) by fixing  two vertices in the graph; one on
each side of the cut. This fixing results in several smaller subproblems that
need to be solved to obtain the new bound. In order to efficiently solve these
subproblems we exploit symmetry in the data; that is, both symmetry in the
min-cut problem and symmetry in the graphs. To obtain upper
bounds for the bandwidth of graphs with symmetry, we develop a heuristic
approach based on the well-known reverse Cuthill-McKee algorithm, and that
improves significantly its performance on the tested graphs. Our approaches
result in the best known lower and upper bounds for the bandwidth of all graphs
under consideration, i.e., Hamming graphs, 3-dimensional generalized Hamming
graphs, Johnson graphs, and Kneser graphs, with up to 216 vertices.
\end{abstract}


\section{Introduction}

For (undirected) graphs, the bandwidth problem (BP) is the problem of labeling the vertices of a given graph with
distinct integers such that the maximum difference between the labels of adjacent vertices is minimal. Determining the
bandwidth is NP-hard (see \cite{Papadem:76}) and it remains NP-hard even if it is restricted to trees with maximum
degree three (see \cite{GagrJoKn:78}) or to caterpillars with hair length three (see \cite{Monien}).

The bandwidth problem originated in the 1950s from sparse matrix computations, and received much attention since
Harary's \cite{Harary:67} description of the problem and Harper's paper \cite{Harper:64} on the bandwidth of the
hypercube. The bandwidth problem arises in many different engineering applications that try to achieve efficient
storage and processing. It also plays a role in designing parallel computation networks, VLSI layout, constraint
satisfaction problems, etc., see, e.g.,~\cite{ChChDeGi:82,DiPeSe:02,LaWi:99}, and the references therein. Berger-Wolf
and Reingold \cite{BergWoRe:02} showed that the problem of designing a code to minimize distortion in multi-channel
transmission can be formulated as the bandwidth problem for the (generalized) Hamming graphs.

The bandwidth problem has been solved for a few families of graphs having special properties. Among these are the path,
the complete graph, the complete bipartite graph \cite{Chvatal1}, the hypercube graph \cite{Harper66}, the grid graph
\cite{Chvalatalova}, the complete $k$-level $t$-ary tree \cite{Smith}, the triangular graph \cite{HwLa77}, and the
triangulated triangle \cite{HoMcDSa}. Still, for many other interesting families of graphs, in particular the
(generalized) Hamming graphs, the bandwidth is unknown. Harper \cite{Harper03} and Berger-Wolf and Reingold
\cite{BergWoRe:02} obtained general bounds for the Hamming graphs, but these bounds turn out to be very weak for
specific examples, as our numerical results will show.

The following lower bounding approaches were recently considered.
Helmberg et al.~\cite{HelReMoPo:95} derived a lower bound for the bandwidth of a graph by exploiting spectral
properties of the graph. The same lower bound was derived by Haemers \cite{haemers} by exploiting interlacing of
Laplacian eigenvalues. Povh and Rendl  \cite{PoRe:07} showed that this eigenvalue bound can also be obtained by solving
a semidefinite programming (SDP) relaxation for the minimum cut (MC) problem. They further tightened the derived SDP relaxation
and consequently obtained a stronger lower bound for the bandwidth. Blum et al.~\cite{BlKoRaVe} proposed a SDP
relaxation for the bandwidth problem that was further exploited by Dunagan and Vempala in \cite{DuVe:01} to derive an
$O(\log^3 n \sqrt{\log \log n })$ approximation algorithm (where $n$ is the number of vertices). De Klerk et
al.~\cite{dKlSoNa} proposed two lower bounds for the graph bandwidth based on SDP relaxations of the quadratic
assignment problem (QAP), and exploited symmetry in some of the considered graphs to solve these relaxations. Their
numerical results show that both their bounds dominate the bound of Blum et al.~\cite{BlKoRaVe}, and that in most of
the cases they are stronger than the bound by Povh and Rendl \cite{PoRe:07}. It is important to remark that all of the
above mentioned SDP bounds are computationally very demanding already for relatively small graphs, that is, for graphs
on about $30$ vertices.  Here, we present a SDP-based lower bound for the bandwidth problem
that dominates the above mentioned bounds and that is also suitable for large graphs with symmetry,
and improve the best known upper bounds for all graphs under consideration. \\

\noindent
{\bf Main results and outline}

\vspace{0.2cm}
\noindent
In this paper we derive a new lower bound for the bandwidth problem of a graph that is based on a new SDP relaxation
for the minimum cut problem. Due to the quality of the new bound for the min-cut problem, and the here improved relation between the min-cut
and bandwidth problem  from \cite{PoRe:07},  we derive the best known lower bounds for the bandwidth of all graphs under consideration, i.e., Hamming
graphs, 3-dimensional generalized Hamming graphs, Johnson graphs, and Kneser graphs. The computed lower bounds for
Hamming graphs turn out to be stronger than the corresponding theoretical bounds by Berger-Wolf and Reingold
\cite{BergWoRe:02}, and Harper \cite{Harper03}.

In particular, all new results one can find  in Section \ref{sec:New}, and they are briefly summarized below.
The new relaxation of the min-cut is a strengthened  SDP relaxation for the QAP by Zhao et
al.~\cite{ZhKaReWo:98} (see also Section \ref{sec:MCQAP}). The new relaxation is obtained from the latter
relaxation by fixing two vertices in the graph; one on each side of the cut (see Section \ref{sec:NEW}). In
Section \ref{sect:McGPP} we show that the new relaxation is equivalent to the SDP relaxation of the graph
partition problem (GPP) from \cite{WolkZhao:99} plus two constraints that correspond to fixing two vertices
in the graph. Although fixing vertices is not a new idea, to the best of our knowledge it is the first time
that one fixes two vertices (an edge or nonedge; instead of a single vertex) in the context of the QAP, GPP,
or MC. Here, it is indeed much more natural to fix two vertices, even though this has some more complicated
technical consequences, as we shall see. Our approach results in several smaller subproblems that need to be
solved in order to obtain the new bound for the min-cut problem. The number of subproblems depends on the
number of `types' of edges and non-edges of the graph, and this number is typically small for the graphs that
we consider. In order to solve the SDP subproblems, we exploit the symmetry in the mentioned graphs and
reduce the size of these problems significantly, see Section \ref{sec:sym2}. We are therefore able to compute
lower bounds for the min-cut of graphs with as much as $216$ vertices, in reasonable time. Finally, to obtain
a lower bound for the bandwidth problem of a graph from the lower bound for the min-cut, we use the new
relation between the mentioned problems from Section \ref{sec:McBnd}.

In order to evaluate the lower bounds, we also compute upper bounds for the bandwidth of the above mentioned
graphs by implementing  a heuristic that improves the well-known reverse Cuthill-McKee algorithm \cite{CutMcKee}, see Section \ref{Sect:CMheuristic}.
Consequently, we are able to determine an optimal labeling (and hence the bandwidth) for several graphs under consideration,
thus showing that for some instances, our lower bound is tight. \\

\noindent
The further set-up of the paper is as follows. In Section \ref{sec:prelim} we introduce notation, provide some
definitions, and give some background on symmetry in graphs.
In Section \ref{sec:Old} we review known bounds for the BP.
In Section \ref{sec:New} we present our new results on obtaining lower and upper bounds for the bandwidth of a graph.
Our numerical results are presented in Section \ref{sec:NumR}.

\section{Preliminaries} \label{sec:prelim}

In this section, we  give some definitions, fix some notation and provide basic information on symmetry in graphs.

\subsection{Notation and definitions}

The space of $k\times k$ symmetric matrices is denoted by $\cS_k$ and the space of $k\times k$ symmetric positive
semidefinite  matrices by $\cS^+_k$. We will sometimes also use the notation $X \succeq 0$ instead of $X \in \cS^+_k$,
if the order of the matrix is clear from the context. For two matrices $X,Y\in \R^{n\times n}$, $X\geq Y$ means
$x_{ij}\geq y_{ij}$, for all $i,j$. The group of $n \times n$ permutation matrices is denoted by $\Pi_n$, whereas the
group of permutations of $\{1,2,\dots,n\}$ is denoted by ${\rm Sym}_n$.

For index sets $\alpha,\beta \subset  \{1,\ldots, n\}$, we denote the submatrix that contains the rows of $A$ indexed
by $\alpha$ and the columns indexed by $\beta$ as $A(\alpha,\beta)$. If $\alpha=\beta$, the principal submatrix
$A(\alpha,\alpha)$ of $A$ is abbreviated as $A(\alpha)$. To denote column $i$ of the matrix $X$ we write $X_{:,i}$.

We use $I_n$ to denote the identity matrix of order $n$, and $e_i$ to denote the  $i$-th standard basis vector.
Similarly, $J_n$ and $u _n$ denote the $n \times n$ all-ones matrix and all-ones $n$-vector, respectively. We will omit
subscripts if the order is clear from the context. We set $E_{ij}=e_i e_j^{\mathrm{T}}$.

The `$\diag$' operator maps an $n\times n$ matrix to the $n$-vector given by its diagonal, while the `$\vect$' operator
stacks the columns of a matrix. The adjoint operator of `diag' we denote by `Diag'. The trace operator is denoted by
`tr'. We will frequently use the property that $\tr AB=\tr BA$.\\

For a graph $G=(V,E)$ with $|V|=n$ vertices,
a {\em labeling} of the vertices of $G$ is a bijection  $\phi:V \to \{1,\ldots,n\}$.
The bandwidth of the labeling $\phi$ of $G$ is defined as
\[
\sigma_{\infty}(G,\phi):=\max_{\{i,j\}\in E} | \phi(i) - \phi(j) |.
\]
The {\em bandwidth} $\sigma_{\infty}(G)$ of a graph $G$ is the minimum of the bandwidth of a labeling of $G$ over all labelings, i.e.,
\[
\sigma_{\infty}(G) := \min \left \{ \sigma_{\infty}(G,\phi)~ | ~\phi: V \to \{1,\ldots,n\} \text{ bijective}  \right \}.
\]
The {\em Kronecker product} $A \otimes B$ of matrices $A \in {\R}^{p \times q}$ and $B\in {\R}^{r\times s}$ is
defined as the $pr \times qs$ matrix composed of $pq$ blocks of size $r\times s$, with block $ij$ given by $a_{ij}B$ $(i = 1,\ldots,p;
j = 1,\ldots,q)$. The following properties of the Kronecker product will be used in the paper, see, e.g.,~\cite{Grah:81} (we assume that the dimensions of the
matrices appearing in these identities are such that all expressions are well-defined):
\begin{eqnarray} \label{kronec1}
(A \otimes B)^\mathrm{T} = A^\mathrm{T} \otimes B^\mathrm{T},~~(A\otimes B)(C\otimes D)=AC \otimes BD.
\end{eqnarray}

\subsection{Symmetry in graphs} \label{sec:SyminGr}
An {\em automorphism} of a graph $G=(V,E)$ is a bijection $\pi:V \rightarrow V$ that preserves edges, that is, such
that $\{\pi(x),\pi(y)\} \in E$ if and only if $\{x,y\}\in E$. The set of all automorphisms of $G$ forms a group under
composition; this is called the {\em automorphism group} of $G$. The {\em orbits} of the action of the automorphism
group acting on $V$ partition the vertex set $V$; two vertices are in the same orbit if and only if there is an
automorphism mapping one to the other. The graph $G$ is {\em vertex-transitive} if its automorphism group acts
transitively on vertices, that is, if for every two vertices, there is an automorphism that maps one to the other (and
so there is just one orbit of vertices). Similarly, $G$ is {\em edge-transitive} if its automorphism group acts
transitively on edges. In this paper, we identify the automorphism group of the graph with the automorphism group of
its adjacency matrix. Therefore, if $G$ has adjacency matrix $A$ we will also refer to the automorphism group of the
graph as $\aut(A) :=\{P \in \Pi_n: P^\mathrm{T}AP=A\}$.

As a generalization of the above algebraic symmetry, combinatorial symmetry is captured in the concept of a coherent
configuration as introduced by Higman \cite{Higmandef} (see also \cite{Higman87}); indeed as a generalization of the
orbitals (the orbits of the action on pairs) of a permutation group. It is defined as follows.

\begin{deff}[Coherent configuration]
\label{def:coherent config} A set of zero-one $n\times n$ matrices $ \{A_1,\ldots, A_r\}$ is called a \emph{coherent
configuration} of rank $r$ if it satisfies the following properties:
\begin{enumerate}[(i)]
\item $\sum_{i \in \mathcal{I}} A_i = I$ for some index set $\mathcal{I} \subset \{1,\ldots,r\}$ and
    $\sum_{i=1}^r A_i = J$,
\item $A_i^\mathrm{T} \in \{A_1,\ldots, A_r\}$ for $i=1,\ldots, r$,
\item $A_iA_j \in \spen\{A_1,\ldots, A_r\}$ for all $i,j$.
\end{enumerate}
\end{deff}
We call $\mathcal{A}:=\spen\{A_1,\dots,A_r\}$ the associated {\em coherent algebra}, and note that this is a
matrix $*$-algebra. If the coherent configuration is commutative, that is, $A_iA_j=A_jA_i$
for all $i,j=1,\dots,r$, then we call it a (commutative) {\it association scheme}. In this case, $\mathcal{I}$ contains
only one index, and it is common to call this index $0$ (so $A_0=I$), and $d:=r-1$ the number of classes of the
association scheme.

One should think of the (nondiagonal) matrices $A_i$ of a coherent configuration as the adjacency matrices of (possibly
directed) graphs on $n$ vertices. The diagonal matrices represent the different `kinds' of vertices (so there are
$|\mathcal{I}|$ kinds of vertices; these generalize the orbits of vertices).  In order to identify the combinatorial
symmetry in a graph, one has to find a coherent configuration (preferably of smallest rank) such that the adjacency
matrix of the graph is in $\mathcal{A}$, see, e.g.,~\cite{KOP}. In that case the nondiagonal matrices $A_i$ represent the different
`kinds' of edges and nonedges.

Every matrix $*$-algebra has a canonical block-diagonal structure. This is a consequence of the theorem by Wedderburn
\cite{Wedderburn} that states that there is a $*$-isomorphism
\[
\varphi:\cA \longrightarrow \oplus_{i=1}^p  \mathbb{C}^{n_i\times n_i}.
\]
Note that in the case of an association schema, all matrices can be diagonalized simultaneously, and the corresponding
$*$-algebra has a canonical diagonal structure $\oplus_{i=0}^d \mathbb{C}$.

We next provide several examples of coherent configurations and association schemes that are used in the remainder of
the paper.

\begin{example}[The cut graph] \label{ex:Kkl}
Let $m=(m_1,m_2,m_3)$ be such that $m_1+m_2+m_3=n$.
The adjacency matrix of the graph $G_{m_1,m_2,m_3}$, which we will call {\em the cut graph}, is given by
\begin{equation} \label{Bmc}
B=
\left (
\begin{array}{ccc}
0_{m_1\times m_1} & J_{m_1\times m_2} & 0_{m_1\times m_3}\\
J_{m_2\times m_1} & 0_{m_2 \times m_2} & 0_{m_2\times m_3}\\
0_{m_3\times m_1} & 0_{m_3 \times m_2} & 0_{m_3\times m_3}
\end{array}
\right ).
\end{equation}
The cut graph is edge-transitive, and belongs to the coherent algebra spanned by the coherent configuration of rank 12
that consists of the matrices
\[
B_1 =
\begin{pmatrix}
I & 0 & 0 \\
0 & 0 & 0 \\
0 & 0 & 0 \\
\end{pmatrix}, \;
B_2 =
\begin{pmatrix}
J-I & 0 & 0 \\
0 & 0 & 0 \\
0 & 0 & 0 \\
\end{pmatrix}, \;
B_3 =
\begin{pmatrix}
0 & J & 0 \\
0 & 0 & 0 \\
0 & 0 & 0 \\
\end{pmatrix}= B_5^{\mathrm{T}},
\]
\[
B_4 =
\begin{pmatrix}
0 & 0 & J \\
0 & 0 & 0 \\
0 & 0 & 0 \\
\end{pmatrix}= B_9^{\mathrm{T}},\;
B_6 =
\begin{pmatrix}
0 & 0 & 0 \\
0 & I & 0 \\
0 & 0 & 0 \\
\end{pmatrix}, \;
B_7 =
\begin{pmatrix}
0 & 0 & 0 \\
0 & J-I & 0 \\
0 & 0 & 0 \\
\end{pmatrix},
\]
\[
B_8 =
\begin{pmatrix}
0 & 0 & 0 \\
0 & 0 & J \\
0 & 0 & 0 \\
\end{pmatrix}= B_{10}^{\mathrm{T}}, \;
B_{11} =
\begin{pmatrix}
0 & 0 & 0 \\
0 & 0 & 0 \\
0 & 0 & I \\
\end{pmatrix},\;
B_{12} =
\begin{pmatrix}
0 & 0 & 0 \\
0 & 0 & 0 \\
0 & 0 & J-I \\
\end{pmatrix},
\]
where the sizes of the blocks are the same as in \eqref{Bmc}.
The coherent algebra is isomorphic to $\mathbb{C} \oplus \mathbb{C}\oplus \mathbb{C}\oplus
\mathbb{C}^{3 \times 3}$. In Appendix \ref{Appendix*isomorph}, the interested reader can find how the associated
$*$-isomorphism $\varphi$ acts on the matrices $B_i$, for $i=1,\ldots, 12$. Note that the cut graph $G_{m_1,m_1,m_3}$
(that is, $m_1=m_2$) is in the coherent configuration of rank $7$ consisting of the matrices $B_1+B_6$, $B_2+B_7$,
$B_3+B_5$, $B_4+B_8$, $B_9+B_{10}$,  $B_{11}$, and $B_{12}$.
\end{example}

\begin{example}[The Hamming graph] \label{ex:hamming}
The Hamming graph $H(d,q)$ is the Cartesian product of $d$ copies of the complete graph $K_q$. The Hamming graph
$H(d,2)$ is also known as the hypercube (graph) $Q_d$. With vertices represented by $d$-tuples of letters from an
alphabet of size $q$, the adjacency matrices of the corresponding association scheme are defined by the number of
positions in which two $d$-tuples differ. In particular, $(A_i)_{x,y}=1$ if $x$ and $y$ differ in $i$ positions;
and the Hamming graph has adjacency matrix $A_1$.

\end{example}

\begin{example}[The generalized Hamming graph] \label{ex:3dimgenham}
The 3-dimensional generalized Hamming graph $H_{q_1,q_2,q_3}$ is the Cartesian product of $K_{q_1}$, $K_{q_2}$, and
$K_{q_3}$. With $V=Q_1 \times Q_2 \times Q_3$, where $Q_i$ is a set of size $q_i$ ($i=1,2,3$), two triples are adjacent
if they differ in precisely one position. Its adjacency matrix is the sum of three adjacency matrices in the
corresponding 7-class association scheme on $V$, that is, the scheme with the $2^3$ adjacency matrices describing being
the same or different for each of the three coordinates.
\end{example}

\begin{example}[The Johnson and Kneser graph] \label{ex:johnson}
Let $\Omega$ be a fixed set of size $v$ and let $d$
be an integer such that $1\leq d\leq v/2$. The vertices of the Johnson scheme are the subsets of $\Omega$ with size $d$.
The adjacency matrices of the association scheme are defined by the size of the intersection of these
subsets, in particular $(A_i)_{\omega,\omega'}=1$ if the subsets $\omega$ and $\omega'$ intersect in $d-i$ elements, for $i=0,\dots,d$.
The matrix $A_1$ represents the Johnson graph $J(v,d)$ (and $A_i$ represents being at distance
$i$ in $G$). For $d=2$, the graph $G$ is strongly regular and also known as a triangular graph.

The Kneser graph $K(v,d)$ is the graph with adjacency matrix $A_d$, that is, two subsets are adjacent whenever they are
disjoint. The Kneser graph $K(5,2)$ is the well-known Petersen graph.
\end{example}

\section{Old  bounds for the bandwidth problem} \label{sec:Old}

The bandwidth problem can be formulated as a  quadratic assignment problem,
but is also closely related to the minimum cut problem that is a special case of the graph partition problem.
In this section we briefly discuss both approaches and present known bounds for the bandwidth problem of a graph.

\subsection{The bandwidth and the quadratic assignment problem} \label{sec:BPQAP}

In terms of matrices, the bandwidth problem asks for a simultaneous permutation of the rows and columns of the adjacency matrix $A$
of the graph $G$ such that all nonzero entries are as close as possible to the main diagonal.
Therefore, a  `natural' problem formulation is the following one.

Let  $k$ be an integer such that $1\leq k \leq n-2$, and let $B=(b_{ij})$ be the $n\times n$ matrix
defined by
\begin{equation}\label{Bbw}
b_{ij} :=
\left \{
\begin{array}{ll}
1 & ~\mbox{for} ~|i-j|>k \\[0.7ex]
0 & ~\mbox{otherwise}.
\end{array} \right .
\end{equation}
Then, if an optimal value of the quadratic assignment problem
\[%begin{equation}
\min_{X \in \Pi_n}  \tr (X^{\mathrm{T}}AX)B \\[1ex]
\]%end{equation}
is zero, then the bandwidth of $G$ is at most $k$.

The idea of formulating  the bandwidth problem as a QAP  where $B$ is defined as above was already suggested by
Helmberg et al.~\cite{HelReMoPo:95} and further exploited by De Klerk et al.~\cite{dKlSoNa}. In \cite{dKlSoNa}, two
SDP-based bounds for the bandwidth problem were proposed; one related to the SDP relaxation for the QAP by Zhao et al.~\cite{ZhKaReWo:98},
and the other one to the improved SDP relaxation for the QAP that is suitable for vertex-transitive
graphs and that was derived by De Klerk and Sotirov \cite{dKSo:12}. The numerical experiments in \cite{dKlSoNa} show that
the latter lower bound for the bandwidth dominates other known bounds for all tested graphs. Although in \cite{dKlSoNa} the
symmetry in the graphs under consideration was exploited, the strongest suggested relaxation was hard to
solve already for graphs with $32$ vertices. This is due to the fact that in general there is hardly any (algebraic) symmetry in
$B$, that is, the automorphism group of $B$ only has order two.

\subsection{The bandwidth and the minimum cut problem} \label{sec:BPMC}

The bandwidth problem is related to the following graph partition problem.
Let $(S_1, S_2, S_3)$ be a partition of $V$ with $|S_i|=m_i$ for $i=1,2,3$. The {\em minimum cut} (MC) problem is:
\[
({\rm MC}) \quad
\begin{array}{lll}
{\rm OPT}_{\rm MC} :=& \min &  \sum\limits_{i\in S_1,j \in S_2} a_{ij} \\[2ex]
          & {\rm s.t.} & (S_1, S_2, S_3) ~~\mbox{partitions}~~ V \\[1.ex]
          & & |S_i|=m_i, ~~i=1,2,3,
\end{array}
\]
where $A=(a_{ij})$ is the adjacency matrix of $G$. To avoid trivialities, we assume that $m_1 \geq 1$ and $m_2 \geq 1$.
We remark that the min-cut problem is known to be NP-hard \cite{GaJo:79}. \\

Helmberg et al.~\cite{HelReMoPo:95} derived a lower bound for the  min-cut problem using the Laplacian eigenvalues of
the graph. In particular, for $m=(m_1,m_2,m_3)$ this bound has closed form expression
\begin{equation}\label{closedForm}
{\rm OPT}_{\rm eig} = -\frac{1}{2} \mu_2 \lambda_2 - \frac{1}{2} \mu_1 \lambda_n,
\end{equation}
where $\lambda_2$ and $\lambda_n$ denote the second smallest and the largest Laplacian eigenvalue of the graph, respectively,
and $\mu_1$ and $\mu_2$ (with $\mu_1\geq \mu_2$) are given by
\[
\mu_{1,2} = \frac{1}{n} \left ( -m_1 m_2 \pm \sqrt{m_1 m_2 (n-m_1)(n-m_2)} \right ).
\]
Further, Helmberg et al.~\cite{HelReMoPo:95} concluded that if ${\rm OPT}_{\rm eig}> 0$ for some $m = (m_1,m_2,m_3)$ then $\sigma(G)\geq m_3+1$.
We remark that the bound on the bandwidth by Helmberg et al.~\cite{HelReMoPo:95} was also derived by Haemers \cite{haemers}. \\

Povh and Rendl \cite{PoRe:07}  proved that ${\rm OPT}_{\rm eig}$ is the solution of a SDP relaxation
of a certain copositive program. They further improved this SDP relaxation by  adding nonnegativity
constraints to the matrix variable. The resulting relaxation is as follows:
\[
({\rm MC}_{\rm COP}) \quad
\begin{array}{lll}
\min & \frac{1}{2} \tr (D \otimes A)Y & \\[1ex]
{\rm s.t.} & \frac{1}{2} \tr( (E_{ij}+E_{ji})\otimes I_n)Y = m_i \delta_{ij}, &~~ 1\leq i\leq j\leq 3\\[1ex]
& \tr(J_3\otimes E_{ii})Y=1, &~~ 1\leq i \leq n \\[1ex]
& \tr( V_i^{\mathrm{T}} \otimes W_j)Y=m_i, &~~ 1\leq i \leq 3, ~~1\leq j\leq n \\[1ex]
& \frac{1}{2} \tr((E_{ij}+E_{ji}) \otimes J_n)Y=m_im_j, &~~ 1\leq i\leq j \leq 3 \\[1ex]
& Y\geq 0, ~Y \in \cS^+_{3n},&
\end{array}
\]
where $D=E_{12} + E_{21}\in \R^{3\times 3}$,
$V_i=e_iu_3^\mathrm{T}\in \R^{3\times 3}$, $W_j=e_ju_n^\mathrm{T}\in \R^{n\times n}$, $1\leq i \leq 3$, $1\leq j
\leq n$, and $\delta_{ij}$ is the Kronecker delta. We use the abbreviation COP to emphasize that ${\rm MC}_{\rm COP}$  is
obtained from a linear program over the cone of completely positive matrices.

In \cite{PoRe:07} Povh and Rendl prove the following proposition, which generalizes
the fact that if ${\rm OPT}_{\rm MC}>0$ for some $m=(m_1,m_2,m_3)$  then $\sigma_{\infty}(G)\geq m_3+1$.

\begin{prop}{\rm \cite{PoRe:07}} \label{PropPoRe}
Let $G$ be an undirected and unweighted graph, and let $m=(m_1,m_2,m_3)$ be such that
${\rm OPT}_{\rm MC}\geq \alpha >0$. Then
\[
\sigma_{\infty}(G) \geq \max \{ m_3+1, m_3+ \lceil \sqrt{2\alpha} \rceil -1 \}.
\]
\end{prop}
\noindent
In the following section we strengthen this inequality.

\section{New lower and upper bounds for the bandwidth problem} \label{sec:New}
In this section we present a new lower bound for the MC problem
that is obtained by strengthening the  SDP relaxation for the GPP or QAP by fixing two vertices in the graph,
and prove that it dominates  ${\rm MC}_{\rm COP}$.
Further, we show how to exploit symmetry in graphs in order to efficiently compute the new MC bound for larger graphs.
We also strengthen Proposition \ref{PropPoRe} that relates the MC and the bandwidth problem,
and present our heuristic that improves the reverse Cuthill-McKee heuristic for the bandwidth problem.

\subsection{The minimum cut and the quadratic assignment problem} \label{sec:MCQAP}

As noted already by Helmberg et al.~\cite{HelReMoPo:95}, the min-cut problem is a special case of the following QAP:
\begin{equation}\label{MCQAP}
\min\limits_{X\in \Pi_n}  \frac{1}{2} \tr X^{\mathrm{T}}AXB,
\end{equation}
where  $A$ is the adjacency matrix of the graph $G$ under consideration,
and $B$ (see \eqref{Bmc}) is the adjacency matrix of the cut graph $G_{m_1,m_2,m_3}$.
Therefore, the following SDP relaxation of the QAP (see \cite{{ZhKaReWo:98},{PoRe:09}}) is also a relaxation for the MC: \label{QAPDSP}
\[
(\mbox{MC}_{\rm QAP})~~~~
\begin{array}{rcl}
\min && \frac{1}{2} {\rm tr}(B\otimes A)Y\\[1ex]
{\rm s.t.} && {\rm tr}(I_n\otimes E_{jj})Y=1, ~~{\rm tr}(E_{jj}\otimes I_n)Y=1, \quad  j=1,\ldots,n\\[1ex]
&& {\rm tr} (I_n\otimes(J_n-I_n)+(J_n-I_n)\otimes I_n)Y=0\\[1ex]
&& {\rm tr}(JY)=n^2 \\[1ex]
&& Y\geq 0, ~~Y \in {\mathcal S}^+_{n^2}.
\end{array}
\]
Note that this is the first time that one uses the above relaxation as a relaxation for the minimum cut problem.
One may easily verify that $\mbox{MC}_{\rm QAP}$ is indeed a relaxation of the QAP
by noting that $Y := \vect(X)\vect(X)^{\mathrm{T}}$ is a feasible point of $\mbox{MC}_{\rm QAP}$ for $X \in \Pi_n$,
and that the objective value of $\mbox{MC}_{\rm QAP}$ at this point $Y$ is precisely $\tr X^{\mathrm{T}}AXB$.
Indeed, the (implicit) assignment constraints $Xu_n=X^{\mathrm{T}}u_n=u_n$ on $X \in \Pi_n$ imply the constraints on
$Y= \vect(X)\vect(X)^{\mathrm{T}}$ involving $E_{jj}$; the sparsity constraints,
i.e., $ {\rm tr} (I_n\otimes(J_n-I_n)+(J_n-I_n)\otimes I_n)Y=0$ follow from the orthogonality conditions $XX^{\mathrm{T}}=X^{\mathrm{T}}X=I_n$; and
the constraint ${\rm tr}(JY)=n^2$ follows from the fact that there are $n$ nonzero elements in the corresponding permutation matrix $X$.

The  matrix $B$, see \eqref{Bmc}, has automorphism group of order $m_1!m_2!m_3!$ when $m_1\neq m_2$ and of order $2(m_1!)^2m_3!$ when
$m_1=m_2$. Since the automorphism group of $B$ is large,  one can exploit the symmetry of $B$ to reduce the size of
$\mbox{MC}_{\rm QAP}$  significantly, see, e.g., \cite{{KOP},{deKlSot:10},{dKSo:12},{dKlSoNa}}.
Consequently, the SDP relaxation $\mbox{MC}_{\rm QAP}$ can be solved much more efficiently when $B$ is defined as in \eqref{Bmc}
than when $B$ is defined as in \eqref{Bbw} (in general, the latter $B$ has only one nontrivial automorphism).

We will next show that   ${\rm MC}_{\rm QAP}$  dominates  ${\rm MC}_{\rm COP}$.
In order to do so, we use the following lemma  from \cite{PoRe:09} that
gives an explicit description of the feasible set of ${\rm MC}_{\rm QAP}$.
\begin{lem}{{\rm \cite[Lemma 6]{PoRe:09}}}\label{validineq}
A matrix
\begin{equation} \label{block form11}
Y =
\left(
\begin{array}{ccc}
Y^{(11)} & \cdots & Y^{(1n)} \\
\vdots &\ddots & \vdots \\
Y^{(n1)} &\cdots& Y^{(nn)}
\end{array}\right) \in \cS_{n^2}^+, \quad Y^{(ij)}\in \R^{n\times n},~~ i,j=1,\ldots,n,
\end{equation}
is feasible for ${\rm MC}_{\rm QAP}$ if and only if $Y$ satisfies
\begin{enumerate}[{\em (i)}]
\item ${\rm tr} (I_n\otimes(J_n-I_n)+(J_n-I_n)\otimes I_n)Y=0$,
\item $\tr Y^{(ii)}=1$ for $1\leq i \leq n$, and $\sum_{i=1}^n \diag(Y^{(ii)})=u$,
\item $u^{\mathrm{T}}Y^{(ij)}= \diag(Y^{(jj)})^{\mathrm{T}}$ for $1\leq i,j \leq n$, and
\item $\sum_{i=1}^n Y^{(ij)}=u\diag(Y^{(jj)})^{\mathrm{T}}$ for $1\leq j \leq n$.
\end{enumerate}
\end{lem}
Now we can prove the following theorem.

\begin{thm} \label{thm:QAPvEig}
Let $G$ be an undirected graph with $n$ vertices and adjacency matrix $A$, and $m_1,m_2,m_3>0$, $m_1+m_2+m_3=n$.
Then the SDP relaxation ${\rm MC}_{\rm QAP}$ dominates the SDP relaxation ${\rm MC}_{\rm COP}$.
\end{thm}

\noindent
\Proof
Let $Y \in \cS_{n^2}^+$ be feasible for $\mbox{MC}_{\rm QAP}$ with block form  (\ref{block form11}).
From $Y$ we construct a feasible point $Z \in \cS_{3n}^+$ for  $\mbox{MC}_{\rm COP}$ in the following way.
First, define blocks
\[
\begin{array}{ccc}
Z^{(11)}=\sum\limits_{i,j=1}^{m_1}Y^{(ij)},
~~Z^{(12)}=\sum\limits_{i=1}^{m_1}\sum\limits_{j=m_1+1}^{m_1+m_2} Y^{(ij)},
~~Z^{(13)}=\sum\limits_{i=1}^{m_1}\sum\limits_{j=m_1+m_2+1}^{n} Y^{(ij)}, \\[3ex]
Z^{(22)}=\sum\limits_{i,j=m_1+1}^{m_1+m_2}Y^{(ij)},
~~Z^{(23)}=\sum\limits_{i=m_1+1}^{m_1+m_2}\sum\limits_{j=m_1+m_2+1}^{n} Y^{(ij)},
~~Z^{(33)}=\sum\limits_{i,j=m_1+m_2+1}^{n}Y^{(ij)},
\end{array}
\]
and then collect these blocks in the matrix
\begin{equation*}
\label{YZW}
Z = \left(
\begin{array}{ccc}
 Z^{(11)} & Z^{(12)} & Z^{(13)} \\
Z^{(21)}  & Z^{(22)} & Z^{(23)} \\
Z^{(31)}  & Z^{(32)}  & Z^{(33)} \\
\end{array} \right),
\end{equation*}
where $Z^{(ji)}= (Z^{(ij)})^{\mathrm{T}} $ for $i<j$.

To prove that $\frac{1}{2} \tr( (E_{ij}+E_{ji})\otimes I_n)Z = m_i \delta_{ij}$ for $1\leq i\leq j\leq 3$, we
distinguish the cases $i=j$ and $i \neq j$. In the first case we have that
\[
\tr( E_{ii}\otimes I_n)Z=\tr Z^{(ii)}=\sum\limits_{i=1}^{m_i} \tr Y^{(ii)}=m_i,
\]
where the last equality follows from Lemma \ref{validineq}
(ii). In the other  case, we have that
\[
\tr( (E_{ij}+E_{ji})\otimes I_n)Z = \tr(Z^{(ij)} + Z^{(ji)} )=0,
\]
where the last equality follows from Lemma \ref{validineq} (i).

The constraint
$\tr(J_3\otimes E_{ii})Z=1$ for $1\leq i \leq n$ follows from  the constraint $\tr (I_n\otimes E_{ii})Y=1$ for ${\rm MC}_{\rm QAP}$ and
the sparsity constraint.
To show that $\tr( V_i^{\mathrm{T}} \otimes W_j)Z=m_i$ for $1\leq i \leq 3$, $1\leq j\leq n$,
we will use Lemma \ref{validineq} (iii) and $\tr(I_n\otimes E_{jj})Y=1$.
In particular, let us assume without loss of generality that $i=1$ and $j=2$. Then
\[
\tr( V_1^{\mathrm{T}} \otimes W_2)Z=
\sum\limits_{i=1}^{m_1} u^{\mathrm{T}} \left ( \sum\limits_{j=1}^{n} Y^{(ij)}_{:,2} \right )
=\sum\limits_{i=1}^{m_1} \left ( \sum\limits_{j=1}^{n} Y^{(jj)}_{2,2} \right )=m_1.
\]
From $u^{\mathrm{T}} Y^{(ij)}u=1$ (see Lemma \ref{validineq} (iii)) it follows that
$\tr((E_{ij}+E_{ji}) \otimes J_n)Z= 2m_im_j$ for $1\leq i\leq j \leq 3$.

It remains to prove that $Z\succeq 0$. Indeed, for every $x\in \R^{3n}$, let $\tilde{x}\in \R^{n^2}$ be defined by
\[
\tilde{x}^{\mathrm{T}} :=
\left ( u_{m_1}^{\mathrm{T}}\otimes x^{\mathrm{T}}_{1:n}, u_{m_2}^{\mathrm{T}}\otimes x^{\mathrm{T}}_{n+1:2n}, u_{m_3}^{\mathrm{T}}\otimes x^{\mathrm{T}}_{2n+1:3n} \right ).
\]
Then
$x^{\mathrm{T}} Z x = \tilde{x}^{\mathrm{T}}Y\tilde{x} \geq 0,$
since  $Y\succeq 0$.
Finally, it follows by direct verification that the objective values coincide for every pair of feasible solutions $(Y,Z)$ that are related as described. \hfill\qed\\

Although our numerical experiments show  that the relaxations ${\rm MC}_{\rm QAP}$ and ${\rm MC}_{\rm COP}$ provide the
same bounds for all test instances, we could not prove that they are equivalent. We remark that we computed
${\rm MC}_{\rm COP}$ only for graphs with at most $32$ vertices (see Section \ref{sec:NumR}), since the computations are very expensive for larger graphs.

\subsection{A new MC relaxation by fixing an edge} \label{sec:NEW}

In this section we strengthen the SDP relaxation $\mbox{MC}_{\rm QAP}$ by adding two constraints that correspond to
fixing two entries $1$ in the permutation matrix of the QAP (\ref{MCQAP}). In other words, the additional constraints
correspond to fixing an (arbitrary) edge in the cut graph, and an edge or a nonedge in the graph $G$. In order to
determine which edge or nonedge in $G$ should be fixed, we consider the action of the automorphism group of $G$ on the
set of ordered pairs of vertices. The orbits of this action are the so-called orbitals, and they represent the
`different' kinds of pairs  of vertices; (ordered) edges, and (ordered) nonedges in $G$ (see also Section \ref{sec:SyminGr}).
Let us assume that there are $t$ such orbitals $\mathcal{O}_h$ ($h=1,2,\dots,t$) of edges and nonedges. We will show that in order to
obtain a lower bound for the original problem, it suffices to compute $t$ subproblems of smaller size. This works
particularly well for highly symmetric graphs, because for such graphs $t$ is relatively small. We formally state
the above idea in the following theorem.

\begin{thm} \label{thmFixnew}
Let $G$ be an undirected graph on $n$ vertices, with adjacency matrix $A$, and $t$ orbitals $\mathcal{O}_h$
$(h=1,2,\dots,t)$ of edges and nonedges. Let $m=(m_1,m_2,m_3)$ be such that $m_1+m_2+m_3=n$. Let $(s_1,s_2)$ be an
arbitrary edge in the cut graph $G_{m_1,m_2,m_3}$ (with adjacency matrix $B$ as defined in \eqref{Bmc}), and
$(r_{h1},r_{h2})$ be an arbitrary pair of vertices in $\mathcal{O}_h$ $(h=1,2,\dots,t)$. Let $\Pi_n(h)$ be the set of
matrices $X \in \Pi_n$ such that $X_{r_{h1},s_1}=1$ and $X_{r_{h2},s_2}=1$ $(h=1,2,\dots,t)$. Then
\[
\min_{X \in \Pi_n} \tr X^{\mathrm{T}}AXB =
\min_{h=1,2,\dots,t} \min_{X \in \Pi_{n}(h)} \tr X^{\mathrm{T}}AXB.
\]
\end{thm}

\proof
Consider the QAP in its combinatorial formulation
\[
\min_{\pi \in {\rm Sym}_n} \sum_{i,j=1}^n a_{\pi(i)\pi(j)}b_{ij}.
\]
From this formulation it is clear that if $\pi$ is an optimal permutation, then for every $\sigma \in {\rm aut}(A)$
also $\sigma \pi$ is optimal. Now let $\pi$ indeed be optimal, let $h$ be such that $(\pi(s_1),\pi(s_2)) \in
\mathcal{O}_h$, and let $\sigma$ be an automorphism of $G$ that maps $(\pi(s_1),\pi(s_2))$ to $(r_{h1},r_{h2})$. Then
$\pi^*:=\sigma \pi$ is an optimal permutation that maps $(s_1,s_2)$ to $(r_{h1},r_{h2})$ and that has objective
\[
\sum_{i,j=1}^n a_{\pi^*(i)\pi^*(j)}b_{ij} =  \tr X^{\mathrm{T}}AXB
\]
for a certain $X \in \Pi_{n}(h)$. This shows that
$$\min_{X \in \Pi_n} \tr X^{\mathrm{T}}AXB \geq
\min_{h=1,2,\dots,t} \min_{X \in \Pi_{n}(h)} \tr X^{\mathrm{T}}AXB,$$
and because the opposite inequality clearly holds, this shows the claimed result.
\hfill\qed

\begin{rem}
In \cite{dKSo:12}, it is shown that in a QAP with automorphism group of $A$ or $B$ acting transitively (that is, at
least one of the corresponding graphs is vertex-transitive), one obtains a global lower bound for the original problem
by fixing {\em one} (arbitrary) entry $1$ in the permutation matrix $X$.
\end{rem}

In the following lemma, we show that when we fix {\em two} entries $1$ in the permutation matrix $X$, we again obtain a QAP,
but that is smaller than the original one.
\begin{lem} \label{branch2}
Let   $X\in \Pi_n$ and $r_1, r_2,s_1,s_2\in \{1,2,\ldots, n\}$ be such that $s_1 \neq s_2$, $X_{r_1,s_1}=1$, and $X_{r_2,s_2}=1$.
Let $\alpha= \{ 1,\ldots, n \}  \backslash \{ r_1, r_2 \}$ and $\beta =\{ 1,\ldots, n \}  \backslash \{ s_1,s_2 \}$, and let $A$ and $B$ be symmetric.
Then
\[
\tr X^{\mathrm{T}}AXB = \tr X(\alpha,\beta)^{\mathrm{T}}(A(\alpha) X(\alpha,\beta)B(\beta) + \hat{C}(\alpha,\beta)) + d,
\]
where
\begin{equation}\label{Chat}
\hat{C}(\alpha,\beta) = 2 A(\alpha, r_1)B(s_1,\beta) + 2 A(\alpha, r_2)B(s_2,\beta),
\end{equation}
and $d = a_{r_1r_1} b_{s_1s_1} + a_{r_2r_2} b_{s_2s_2} + 2a_{r_1r_2} b_{s_1s_2}$.
\end{lem}

\proof We will show this by using the combinatorial formulation of the QAP, splitting its summation appropriately
(while using that $A$ and $B$ are symmetric matrices), and then switching back to the trace formulation (we will omit details).
Indeed, let $\pi$ be the permutation that corresponds to $X$; in particular we have that $\pi(s_1)=r_1$ and
$\pi(s_2)=r_2$. Then
\begin{align*}\pushQED{\qed}
\tr X^{\mathrm{T}}AXB &= \sum_{i,j=1}^n a_{\pi(i)\pi(j)}b_{ij} \\
&= \sum_{\{i,j\}\neq \{s_1,s_2\}} a_{\pi(i)\pi(j)}b_{ij} + 2\sum_{j \neq s_1,s_2} a_{r_1\pi(j)}b_{s_1j} + 2\sum_{j \neq s_1,s_2} a_{r_2\pi(j)}b_{s_2j}+d\\
&= \tr X(\alpha,\beta)^{\mathrm{T}}A(\alpha) X(\alpha,\beta)B(\beta)\\
& + 2 \tr X(\alpha,\beta)^{\mathrm{T}}(A(\alpha, r_1)B(s_1,\beta) +  A(\alpha, r_2)B(s_2,\beta))+d\\
&=\tr X(\alpha,\beta)^{\mathrm{T}}(A(\alpha) X(\alpha,\beta)B(\beta) + \hat{C}(\alpha,\beta)) + d.\qedhere \popQED
\end{align*}
Since  $A(\alpha), B(\beta)\in \cS_{n-2}$, where $\alpha$, $\beta$ are defined as in Lemma \ref{branch2},
and $X(\alpha,\beta)\in \Pi_{n-2}$, the reduced problem
\begin{equation}\label{redQAP2}
\min\limits_{ X\in \Pi_{n-2} }\tr X^{\mathrm{T}}(A(\alpha) X B(\beta) + \hat{C}(\alpha,\beta))
\end{equation}
is also a quadratic assignment problem.
Therefore computing the new  lower bound for the min-cut problem
using  Theorem \ref{thmFixnew} reduces to solving several SDP subproblems of the form \label{Mcfixhh}
\[
(\mbox{MC}_{\rm fix}^h)~~~~
\begin{array}{rl}
\mu^*_h=\min & \frac{1}{2} {\rm tr}( (B(\beta)\otimes A(\alpha^h)) + \Diag(\hat{c}))Y + \frac{1}{2}d^h\\[1ex]
{\rm s.t.} & {\rm tr}(I\otimes E_{jj})Y=1, ~~{\rm tr}(E_{jj}\otimes I)Y=1, \quad  j=1,\ldots,n-2\\[1ex]
& {\rm tr} (I\otimes(J-I))+(J-I)\otimes I)Y=0\\[1ex]
& {\rm tr}(JY)=(n-2)^2 \\[1ex]
& Y\geq 0, ~~Y\in \cS^+_{(n-2)^2},
\end{array}
\]
where $h =1,\ldots, t$,  $I,J,E_{jj} \in {\mathcal S}_{n-2}$, $\hat{c}={\rm vec}(\hat{C}(\alpha,\beta))$,
$\beta =\{ 1,\ldots, n \}  \backslash \{ s_1,s_2 \}$,
 $\alpha^h= \{ 1,\ldots, n \}  \backslash \{ r_{h1}, r_{h2} \}$,
and the constant $d^h$ is defined as $d$ in Lemma \ref{branch2}.
Finally, the new lower bound for the min-cut problem is
\[
\mbox{MC}_{\rm fix} = \min_{h=1,\ldots, t }  \mu^*_h.
\]
Thus, computing the new lower bound for the bandwidth problem   $\mbox{MC}_{\rm fix}$ involves solving $t$ subproblems.
This seems like  a computationally  demanding and very restrictive approach, knowing that in general it is hard to
solve $\mbox{MC}_{\rm fix}^h$ already when $n$ is about 15, see \cite{ReSo}.
However, our aim here is to compute lower bounds on the bandwidth problem for graphs that
are  known to be highly symmetric. Consequently, for such graphs $t$ is small and we can exploit the symmetry of graphs
to reduce the size of $\mbox{MC}_{\rm fix}^h$  as described in Section \ref{sec:sym2}.

It remains here to show that the new bound $\mbox{MC}_{\rm fix}$ dominates all other mentioned lower bounds for the min-cut problem.
To show this we need the following proposition.
Note that we already know that the eigenvalue bound  $\mbox{OPT}_{\rm eig}$ is dominated by $\mbox{MC}_{\rm COP}$ (see \cite{PoRe:07}),
which in turn is dominated by $\mbox{MC}_{\rm QAP}$ (see Theorem  \ref{thm:QAPvEig}).
\begin{prop} \label{equivPr}
Let $(s_1,s_2)$ be an arbitrary edge in the cut graph $G_{m_1,m_2,m_3}$
and $(r_{h1},r_{h2})$ be an arbitrary pair of vertices in $\mathcal{O}_h$  $(h=1,2,\dots,t)$.
Also, let
$\alpha^h= \{ 1,\ldots, n \}  \backslash \{ r_{h1}, r_{h2} \}$
and
$\beta =\{ 1,\ldots, n \}  \backslash \{ s_{1},s_{2} \}$.  Then the semidefinite program
\[
({\rm MC}_{\rm {QAP}}^h)~~~~
\begin{array}{rl}
\min & \frac{1}{2} {\rm tr}(B\otimes A)Y\\[1ex]
{\rm s.t.} & {\rm tr}(I_n\otimes E_{jj})Y=1, ~~{\rm tr}(E_{jj}\otimes I_n)Y=1, \quad  j=1,\ldots,n\\[1ex]
& \tr (I_n\otimes(J_n-I_n))+(J_n-I_n)\otimes I_n)Y=0\\[1ex]
& \tr(JY)=n^2 \\[1ex]
& \tr(E_{s_is_i}\otimes E_{r_{hi},r_{hi}})Y=1, \qquad i=1,2              \\[1ex]
&  Y\geq 0, ~~Y \in \cS^+_{n^2}
\end{array}
\]
is equivalent to {\rm MC}$_{\rm fix}^h$  $(h=1,2,\dots,t)$
in the sense that there is a bijection between the feasible sets that preserves the objective function.
\end{prop}

\proof The proof is similar to the proof of Theorem 27.3 in \cite{Sot10}. \hfill\qed\\

\noindent
Now, from Proposition \ref{equivPr} it follows indeed that the new bound dominates all others.
\begin{cor}\label{cor:newbounddominates}
Let $(s_1,s_2)$ be an arbitrary edge in the cut graph $G_{m_1,m_2,m_3}$
and $(r_{h1},r_{h2})$ be an arbitrary pair of vertices in $\mathcal{O}_h$  $(h=1,2,\dots,t)$.
Then the SDP relaxation  ${\rm MC}_{\rm fix}$ dominates  ${\rm MC}_{\rm QAP}$.
\end{cor}

In this section, we proved that the proposed new bound for the min-cut dominates  the SDP bound $\mbox{MC}_{\rm COP}$,
which in turn dominates $\mbox{OPT}_{\rm eig}$.
Numerical experiments by Povh and Rendl \cite{PoRe:07} and De Klerk et al.~\cite{dKlSoNa} show that it is already hard to solve   $\mbox{MC}_{\rm COP}$ for
graphs whose size is $32$. In the following section we show how one can  exploit symmetry of the considered graphs
to efficiently compute the new relaxation $\mbox{MC}_{\rm fix}$ for large instances.

\subsection{Reduction by using symmetry} \label{sec:sym2}

Computational  experiments show that in general, the SDP relaxation of the QAP by Zhao et al.~\cite{ZhKaReWo:98},
i.e., $\mbox{MC}_{\rm QAP}$, cannot be solved in a straightforward
way by interior point methods for instances where $n$ is larger than 15,  see, e.g., \cite{ReSo}.
However, it is possible to solve larger  problem instances when the data matrices  have
large automorphism groups, as described in \cite{{deKlSot:10},{KOP}}.
Further, in \cite{{dKSo:12},{dKPaDoSo:10},{dKlSoNa}} it is shown how to reduce the size of
the SDP relaxation that is obtained from the relaxation from \cite{ZhKaReWo:98} after fixing one entry $1$ in the permutation matrix.
Here, we go one step further and show how to exploit symmetry in the data to
solve efficiently the SDP relaxation that is obtained from $\mbox{MC}_{\rm QAP}$  after fixing two entries $1$ in the
permutation matrix, i.e., $\mbox{MC}_{\rm fix}^h$.

In \cite{{deKlSot:10},{KOP}} it is shown  that if one (or both) of the data matrices belong to a matrix  $*$-algebra, then one can
exploit the structure of the algebra  to reduce the size of the SDP relaxation by Zhao et al.~\cite{ZhKaReWo:98}.
In particular, since the matrix  $B$  belongs to the coherent algebra described in Example \ref{ex:Kkl},
$\mbox{MC}_{\rm QAP}$ reduces to the following:
\begin{equation} \label{reducedSDP}
\begin{array}{rl}
\min \ & \frac{1}{2} \tr A (X_3 + X_5) \\ [1ex]
{\rm s.t.} \ &  X_1 + X_6 + X_{11} =I_n\\ [1ex]
& \sum\limits_{i=1}^{12} X_i=J_n\\     [2ex]
& \tr (J X_i)=p_i,~X_i\ge 0,\quad i=1,\dots,12\\   [1.5ex]
& X_1 - \frac{1}{m_1-1} X_2 \succeq 0, ~X_6 - \frac{1}{m_2-1} X_7 \succeq 0,~X_{11} - \frac{1}{m_3-1} X_{12} \succeq 0\\[2.5ex]
& \begin{pmatrix}
\frac{1}{m_1} (X_1+X_2) & \frac{1}{\sqrt{m_1m_2}} X_3 & \frac{1}{\sqrt{m_1m_3}} X_4\\[1.ex]
\frac{1}{\sqrt{m_1m_2}} X_5 & \frac{1}{m_2} (X_6+X_7) &  \frac{1}{\sqrt{m_2m_3}} X_8 \\[1.ex]
\frac{1}{\sqrt{m_1m_3}} X_9 & \frac{1}{\sqrt{m_2m_3}} X_{10} & \frac{1}{m_3} (X_{11}+X_{12})
\end{pmatrix}  \succeq 0\\ [5.ex]
& X_3=X_5^{\mathrm{T}}, X_4=X_9^{\mathrm{T}}, X_8=X_{10}^{\mathrm{T}},\\ [1.5ex]
&  ~X_1,X_2,X_6,X_7,X_{11},X_{12}\in \cS_{n},
\end{array}
\end{equation}
where $p_1=m_1$, $p_2=m_1(m_1-1)$, $p_3=p_5=m_1 m_2$, $p_4=p_9=m_1 m_3$,  $p_6=m_2$, $p_7=m_2(m_2-1)$,
$p_8=p_{10}=m_2m_3$, $p_{11}=m_3$, and $p_{12}=m_3(m_3-1)$. We used a reformulation of $\mbox{MC}_{\rm QAP}$
from \cite[p.~186]{KOP} and the coherent algebra from Example \ref{ex:Kkl}, and then applied the associated
$*$-isomorphism given in Appendix \ref{Appendix*isomorph} to derive the four SDP constraints in
\eqref{reducedSDP}.

It is possible to perform further reduction of \eqref{reducedSDP} for graphs with symmetry, see, e.g., \cite{dKSo:12}.
In order to do so, we need a coherent algebra containing $A$.
For the graphs that we consider, such coherent algebras are known,  see Examples \ref{ex:hamming}--\ref{ex:johnson}.
Therefore, further reduction of (\ref{reducedSDP}) follows  directly from results from, e.g., \cite{deKlSot:10,dKSo:12}.

The symmetry reduction of  the SDP relaxation ${\rm MC}_{\rm fix}^h$ ($h=1,\ldots,t$)
is however not a straightforward application of results from \cite{{deKlSot:10},{dKSo:12},{dKPaDoSo:10},{dKlSoNa}}.
Therefore we provide a detailed analysis below.
In order to perform the desired symmetry reduction of ${\rm MC}_{\rm fix}^h$,
we will find a large enough subgroup of the automorphism group
of the data matrix from the objective function, i.e., of  $B(\beta)\otimes A(\alpha^h) + \Diag({\rm vec}(\hat{C}))$,
where $\hat{C}$ is given in \eqref{Chat}.
Note that if the objective function would not have a linear term, then ${\rm aut}(B(\beta)) \otimes {\rm aut}(A(\alpha^h))$
could serve as such a subgroup.

In what follows, we first determine appropriate subgroups of the automorphism groups of $A(\alpha)$, $B(\beta)$, and  $\hat{C}(\alpha,\beta)$.
To describe these we introduce the following definition.

\begin{deff}
For  $r_1,r_2\in \{1,\ldots,n \}$ with $r_1\neq r_2$, the
subgroup of  $\mbox{\rm aut}(A)$ that fixes row and column $r_1$ and row and column $r_2$ of $A$ is:
\begin{equation} \label{stabil12}
{\rm stab}((r_1,r_2),A):= \{ P \in  {\rm aut}(A): P_{r_1,r_1}=P_{r_2,r_2}=1 \}.
\end{equation}
This is the {\em pointwise} stabilizer subgroup of aut($A$) with respect to $(r_1,r_2)$, see, e.g., \cite{{Cam:99}}.
\end{deff}
It is important to distinguish the pointwise stabilizer from the {\em setwise} stabilizer of a set $S\subset V$; in the
latter it is not required that every point of $S$ is fixed, but that the set $S$ is fixed (that is, $\pi(S)=S$). For
$\alpha= \{ 1,\ldots, n \} \backslash \{ r_1,r_2 \}$ we define
\[
{\mathcal H}(A(\alpha)) :=\{ P(\alpha): P \in  {\rm stab}((r_1,r_2),A) \}.
\]
This group is a subgroup of the automorphism group of $A(\alpha)$.

Similarly, we define ${\mathcal H}(B(\beta)):=\{ P(\beta): P \in  {\rm stab}((s_1,s_2),B) \}$ for $\beta =\{ 1,\ldots,n \}  \backslash$ $\{ s_1,s_2 \}$. \\

\begin{lem}\label{claim}
The action of ${\mathcal H}(B(\beta))$ has  $12$ orbitals.
\end{lem}
\noindent
{\em Proof.}
In the following, we will consider an edge $(s_1,s_2)$ in the cut graph (without loss of generality $s_1\in S_1$ and $s_2\in S_2$).
Because of the simple structure of the cut graph, ${\mathcal H}(B(\beta))$ can be easily described.
Indeed, because $s_1$ and $s_2$ are fixed by $P \in {\rm stab}((s_1,s_2),B)$, the sets
$S_1\setminus\{s_1\}$ and $S_2\setminus\{s_2\}$ are fixed (as sets) by $P(\beta) \in {\mathcal H}(B(\beta))$. This
implies that ${\mathcal H}(B(\beta))$ is the direct product of the symmetric groups on $S_1\setminus\{s_1\}$,
$S_2\setminus\{s_2\}$, and $S_3$. In fact, this is the full automorphism group of $B(\beta)$ in case $m_1\neq m_2$ (in
the case that $m_1=m_2$ it is an index 2 subgroup of ${\rm aut}(B(\beta))$ since the `swapping' of
$S_1\setminus\{s_1\}$ and $S_2\setminus\{s_2\}$ is not allowed).
Therefore, the action of ${\mathcal H}(B(\beta))$ on $\beta$ has 12 orbitals, similar as described in Example \ref{ex:Kkl}. \hfill\qed \\

Now we can describe the group that we will exploit to reduce the size of the subproblem ${\rm MC}_{\rm fix}^h$ ($h=1,\ldots,t$).

\begin{prop} \label{ThmObjsubGr}
Let $r_1,r_2,s_1,s_2\in \{1,\ldots,n \}$, $\alpha= \{ 1,\ldots, n \}  \backslash \{ r_1, r_2 \}$, $\beta =\{ 1,\ldots,
n \}  \backslash$ $ \{ s_1,s_2 \}$, and $\hat{C}(\alpha,\beta) = 2 A(\alpha, r_1)B(s_1,\beta) + 2 A(\alpha,
r_2)B(s_2,\beta)$. Then
\[
{\mathcal H}(B(\beta)) \otimes {\mathcal H}({A(\alpha)})
\]
is a subgroup of the automorphism group of $B(\beta) \otimes A(\alpha)  + {\rm Diag}({\rm
vec}(\hat{C}(\alpha,\beta)))$.
\end{prop}

\noindent \proof Let $P_B\in {\mathcal H}(B(\beta))$  and $P_A \in {\mathcal H}(A(\alpha))$. It is clear that $P_B
\otimes P_A$ is an automorphism of $B(\beta) \otimes A(\alpha)$, so we may restrict to showing that is also an
automorphism of ${\rm Diag}({\rm vec}(\hat{C}(\alpha,\beta)))$. In order to show this, we will use that for $i=1,2$ we have that
$$P_A^{\mathrm{T}} {\rm Diag}(A(\alpha, r_i)) P_A ={\rm Diag}(A(\alpha, r_i)) \quad {\rm and} \quad
P_B^{\mathrm{T}} {\rm Diag}(B(\beta,s_i)) P_B ={\rm Diag}(B(\beta,s_i)).$$
Indeed, the first equation is equivalent to the (valid) property that $a_{\pi(j)r_i}=a_{jr_i}$ for all $j \neq r_1,r_2$
and all automorphisms $\pi$ of $A$ that fix both $r_1$ and $r_2$, and the second equation is similar. Because $ {\rm
vec}(\hat{C}(\alpha, \beta)) =2B(\beta,s_1) \otimes A(\alpha, r_1)+ 2B(\beta, s_2) \otimes A(\alpha, r_2)$, the result
now follows from
\[
\begin{array}{rl}
& (P_B \otimes P_A)^{\mathrm{T}} {\rm Diag} \left [B(\beta,s_1) \otimes A(\alpha, r_1)
+  B(\beta, s_2) \otimes A(\alpha, r_2) \right ] (P_B \otimes P_A)~~~~~~~~~~~~~~~~~~~~~~~~~~~~~~~~~~~~~~ \\[1.5ex]
=& P_B^{\mathrm{T}} {\rm Diag}(B(\beta,s_1)) P_B \otimes P_A^{\mathrm{T}} {\rm Diag}(A(\alpha, r_1))P_A  \\[1.5ex]
& +~ P_B^{\mathrm{T}} {\rm Diag}(B(\beta, s_2)) P_B \otimes P_A^{\mathrm{T}} {\rm Diag}(A(\alpha, r_2))P_A\\[1.5ex]
=&
{\rm Diag}[B(\beta,s_1)  \otimes A(\alpha, r_1) +  B(\beta, s_2) \otimes A(\alpha, r_2)],
\end{array}
\]
where we have used the properties in (\ref{kronec1}) of the Kronecker product.
\hfill \qed\\

Now, for an associated  matrix $*$-algebra of $A(\alpha)$  one can take
the centralizer ring (or commutant) of ${\mathcal H}({A(\alpha)})$, i.e.,
\[
{\mathcal A}_{A(\alpha)}=\{ X\in \R^{n\times n}: XP=PX, ~\forall P \in {\mathcal H}({A(\alpha)}\}.
\]
Similarly, we take the centralizer ring ${\mathcal A}_{B(\beta)}$ of ${\mathcal H}({B(\beta)})$.
Now we restrict the variable $Y$ in $\mbox{MC}_{\rm fix}^h$ to lie in ${\mathcal A}_{B(\beta)} \otimes {\mathcal A}_{A(\alpha)}$,
and obtain a basis of this algebra from the orbitals of   $ {\mathcal H}({A(\alpha)})$ and $ {\mathcal H}({B(\beta)})$.
Finally, the symmetry reduction of $\mbox{MC}_{\rm fix}^h$ is similar to the
symmetry reduction of $\mbox{MC}_{\rm QAP}$, see \eqref{reducedSDP} and also  \cite{{dKSo:12},{dKPaDoSo:10}}.
The interested reader can find the reduced formulation of $\mbox{MC}_{\rm fix}^h$ in Appendix \ref{AppendixMcfixh}.

\subsection{The minimum cut and the graph partition} \label{sect:McGPP}

In this section we relate ${\rm MC}_{\rm QAP}$ and  ${\rm MC}_{\rm fix}$ with corresponding SDP relaxations for the graph partition problem.
In particular, we show that the strongest  known SDP relaxation for the GPP applied to the MC problem is equivalent to ${\rm MC}_{\rm QAP}$,
and that ${\rm MC}_{\rm fix}$ can be obtained by computing a certain SDP relaxation for the GPP.

As mentioned before, the minimum cut problem is a special case of the graph partition problem.
Therefore one can solve the SDP relaxation for the graph partition problem from Wolkowicz and Zhao \cite{WolkZhao:99}
to obtain a lower bound for the minimum cut.  The relaxation for the GPP from \cite{WolkZhao:99} is as follows:
\[
(\mbox{MC}_{\rm GPP})~~~~
\begin{array}{rl}
\min & \frac{1}{2} \tr (D \otimes A) Y \\[1.5ex]
{\rm s.t.} & \tr( (J_3-I_3)\otimes I_n )Y =0\\
& \tr(I_3 \otimes J_n)Y + \tr(Y) = -(\sum\limits_{i=1}^km_i^2+n)+ 2y^{\mathrm{T}}((m+u_3)\otimes u_{n})\\
& \left ( \begin{array}{cc} 1 & y^{\mathrm{T}} \\ y & Y  \end{array} \right ) \in \cS_{3n+1}^+, \quad Y\geq 0,
\end{array}
\]
where $D=E_{12}+E_{21} \in \R^{3\times 3}$ and $A$ is the adjacency matrix of the graph.
Actually, the SDP relaxation from  \cite{WolkZhao:99} does not include nonnegativity constraints,
but we add them to strengthen the bound, see also \cite{Sot11}.
Note that for the $k$-partition we define $D:=J_k-I_k$.
The results in \cite{Sot11} show that $\mbox{MC}_{\rm GPP}$ is the strongest known SDP relaxation for the graph partition problem.
In the following  theorem we prove that $\mbox{MC}_{\rm GPP}$ is equivalent to  $\mbox{MC}_{\rm QAP}$.
\begin{thm}\label{thm:eqGPPQAP}
Let $G$ be an undirected graph with $n$ vertices and adjacency matrix $A$, and $m=(m_1,m_2,m_3)$ be such that $m_1+m_2+m_3=n$.
Then the SDP relaxations  ${\rm MC}_{\rm GPP}$ and  ${\rm MC}_{\rm QAP}$ are equivalent.
\end{thm}

\noindent
\Proof
Let $Y \in \cS_{n^2}^+$ be feasible for $\mbox{MC}_{\rm QAP}$ with block form (\ref{block form11}) where $Y^{(ij)}\in \R^{n\times n}$, $i,j=1,\ldots,n$.
We construct from $Y\in \cS_{n^2}^+$ a feasible point $(W,w)$ for $\mbox{MC}_{\rm GPP}$ in the following way.
First, define blocks
\begin{equation} \label{r2}
\begin{array}{cc}
W^{(11)}  := \sum\limits_{i,j=1}^{m_1} Y^{(ij)},~~
& W^{(12)}  := \sum\limits_{i=1}^{m_1} \sum\limits_{j=m_1+1}^{m_1+m_2}  Y^{(ij)},\\[3ex]
 W^{(13)}  := \sum\limits_{i=1}^{m_1} \sum\limits_{j=m_1+m_2+1}^{n}  Y^{(ij)},
& W^{(22)}  := \sum\limits_{i,j =m_1+1}^{m_1+m_2}  Y^{(ij)},\\[3ex]
 W^{(23)}  := \sum\limits_{i=m_1+1}^{m_1+m_2} \sum\limits_{j=m_1+m_2+1}^{n}  Y^{(ij)},
& W^{(33)}  := \sum\limits_{i,j=m_1+m_2+1}^{n}   Y^{(ij)},
\end{array}
\end{equation}
and then collect all blocks into the matrix
\begin{equation} \label{YZW}
W := \left (
\begin{array}{ccc}
 W^{(11)} & W^{(12)} & W^{(13)} \\
 (W^{(12)})^{\mathrm{T}} & W^{(22)} & W^{(23)} \\
(W^{(13)})^{\mathrm{T}}  & (W^{(23)})^{\mathrm{T}}  &  W^{(22)}
\end{array} \right ).
\end{equation}
Define $w :=\diag(W)$. It is not hard to verify that $(W,w)$ is feasible for $\mbox{MC}_{\rm GPP}$.
Also, it  is straightforward  to see (by construction) that the two objective values are equal.

Conversely, let  $(W,w)$ be feasible for  $\mbox{MC}_{\rm GPP}$ and suppose that $W$ has the block form (\ref{r2}).
To show that $\mbox{MC}_{\rm GPP}$  dominates  $\mbox{MC}_{\rm QAP}$, we exploit  the fact that  $\mbox{MC}_{\rm QAP}$
reduces to  (\ref{reducedSDP}). This reduction is due to the symmetry  in the corresponding cut graph.
Now, let
\[
\begin{array}{lll}
X_1 := \Diag(\diag(W^{(11)})), ~~& X_2 :=W^{(11)}-X_1, & X_3 := W^{(12)} \\[1ex]
X_4 := W^{(13)}, ~~ & X_6 :=\Diag(\diag(W^{(22)})), ~~& X_7 :=W^{(22)}-X_6\\[1ex]
X_8 := W^{(23)}, ~~& X_{11} := \Diag(\diag(W^{(33)})),  & X_{12} :=W^{(33)} -X_{11}.   \\
\end{array}
\]
It is easy to verify that so defined $X_1$, \ldots, $X_{12}$ are feasible for (\ref{reducedSDP}), thus for $\mbox{MC}_{\rm QAP}$.
It is also easy to see that the two objectives coincide. \hfill\qed\\

The previous theorem proves that  $\mbox{MC}_{\rm GPP}$ and $\mbox{MC}_{\rm QAP}$ provide the same lower bound  for the given minimum cut problem.
Note that the positive semidefinite matrix variable in $\mbox{MC}_{\rm GPP}$  has order $3n+1$, while the largest linear matrix inequality
in \eqref{reducedSDP} has order $3n$. Since the additional row and column in $\mbox{MC}_{\rm GPP}$ makes a symmetry reduction more difficult,  we
choose to compute \eqref{reducedSDP} instead of $\mbox{MC}_{\rm GPP}$.

In the sequel we address the issue of fixing an edge in the graph partition formulation of the MC problem.
In Proposition \ref{equivPr} we proved that $\mbox{MC}_{\rm fix}^h$  can be obtained from  $\mbox{MC}_{\rm QAP}$
by adding the constraints $\tr(E_{s_is_i} \otimes E_{r_{hi}r_{hi}})=1$ for appropriate $(s_1,s_2)$, $(r_{h1},r_{h2})$.
Similarly, we can prove that a SDP relaxation for the GPP formulation of the MC with fixed  $(r_{h1},r_{h2})\in {\mathcal O}_h$ ($h=1,\ldots, t$)
can be obtained from $\mbox{MC}_{\rm GPP}$ by adding the constraints
\[
\tr(\bar{E}_{ii} \otimes E_{r_{hi}r_{hi}})Y=1, \quad i=1,2,
\]
where $\bar{E}_{ii} \in \R^{3\times 3}$ and $E_{r_{hi}r_{hi}}\in \R^{n\times n}$.
We finally arrive at the following important result.
\begin{thm} \label{eq:QAPGPPfix}
Let $(r_{h1},r_{h2})$ be an arbitrary pair of vertices in $\mathcal{O}_h$  $(h=1,2,\dots,t)$.
 Then the semidefinite program
\[
({\rm MC}_{\rm {GPP}}^h)~~~~
\begin{array}{rl}
\min & \frac{1}{2} \tr (D \otimes A) Y \\[1.5ex]
{\rm s.t.} & \tr( (J_3-I_3)\otimes I_n )Y =0\\
& \tr(I_3 \otimes J_n)Y + \tr(Y) = -(\sum\limits_{i=1}^km_i^2+n)+ 2y^{\mathrm{T}}((m+u_3)\otimes u_{n})\\
& \tr(\bar{E}_{ii} \otimes E_{r_{hi}r_{hi}})Y=1, \quad i=1,2 \\[2ex]
& \left ( \begin{array}{cc} 1 & y^{\mathrm{T}} \\ y & Y  \end{array} \right ) \in \cS_{3n+1}^+, \quad Y\geq 0,
\end{array}
\]
is equivalent to {\rm MC}$_{\rm QAP}^h$
in the sense that there is a bijection between the feasible sets that preserves the objective function.
\end{thm}

\noindent
{\em Proof}. The proof is similar to  the proof of Theorem \ref{thm:eqGPPQAP}.
However, one should take into consideration that the constraint $\tr(E_{s_1s_1}\otimes E_{r_{h1}r_{h1}})Y=1$ reduces  to
$\tr (E_{r_{h1}r_{h1}} X_1)=1$ in \eqref{reducedSDP},
and  $\tr(E_{s_2s_2}\otimes E_{r_{h2}r_{h2}})Y=1$ to  $\tr (E_{r_{h2}r_{h2}} X_6)=1$. \hfill\qed\\

This theorem shows that we can obtain MC$_{\rm fix}$ also by computing ${\rm MC}_{\rm {GPP}}^h$ for all $h=1,\ldots,t$.

\subsection{On computing the new lower bound for the bandwidth} \label{sec:McBnd}

In this section we improve the Povh-Rendl inequality from  \cite{PoRe:07}
that relates the minimum cut and the bandwidth problem. In particular, we derive the following result.

\begin{prop}\label{computeBdw}
Let $G$ be an undirected and unweighted graph, and let $m=(m_1,m_2,m_3)$ be such that
${\rm OPT}_{\rm MC}\geq \alpha >0$. Then
\[
\sigma_{\infty}(G) \geq m_3 + \left \lceil -\frac{1}{2}+ \sqrt{2\lceil \alpha \rceil +\frac{1}{4} }\right \rceil.
\]
\end{prop}

\noindent
\Proof  (See also  \cite{PoRe:07}.)  Let $\phi$ be an optimal labeling of $G$, and let $(S_1,S_2,S_3)$ be a partition of $V$,
such that $\phi(S_1)=\{1,\ldots, m_1\}$ and $\phi(S_2)=\{m_1+m_3+1,\ldots, n \}$.
Let $\Delta$ be the maximal difference of labels over all edges connecting sets $S_1$ and $S_2$ (so $\Delta \leq \sigma_{\infty}(G)$)
and let $\delta=\Delta -m_3$ (and note that $\delta \geq 1$).
Since the number of edges between $S_1$ and $S_2$ is  at most $\delta(\delta+1)/2$, it follows that $\delta(\delta+1)\geq 2\cdot {\rm OPT}_{\rm MC}$.
Note first of all that this implies Proposition \ref{PropPoRe}.
Secondly, from  $\delta(\delta+1)\geq 2 \lceil \alpha \rceil$ and $\sigma_{\infty}(G) \geq m_3 + \delta$, the required inequality follows.  \hfill\qed\\

In  Section \ref{sec:NumR}, we compute the new bound for the BP  by using Proposition \ref{computeBdw}.
It is worth mentioning that by using the  expression from Proposition \ref{computeBdw}
instead of the expression from Proposition \ref{PropPoRe}, the bounds may improve by several integer values.
The largest improvement that we recorded is $3$ integer values.

\subsection{An improved reverse Cuthill-McKee algorithm}\label{Sect:CMheuristic}

In order to get more information about the quality of the lower bounds, we needed good upper bounds for the bandwidth of a graph.
We obtained these by testing the well known (reverse) Cuthill-McKee algorithm \cite{CutMcKee} on several graphs with symmetry;
however the output seemed far from optimal. This is not at all surprising because the Cuthill-McKee algorithm sorts mostly
on vertex degrees, and in the graphs of our interest these are all equal.
Therefore we developed a heuristic that combines the reverse Cuthill-McKee algorithm and an improvement procedure.
The details of this improvement procedure are described as follows.

Consider a labeling $\phi$ of the graph with (labeling) bandwidth $\sigma_{\infty}$. A vertex $u$ is called critical if it has a neighbor
$w$ such that $|\phi(u)-\phi(w)|=\sigma_{\infty}$. Now we consider the critical vertex $u$ that has the largest label, and its
critical neighbor $w$. Let $z$ be the vertex with the largest label that is not adjacent to any of the vertices with
labels $1,2,\dots,\phi(w)$. If $\phi(z)<\phi(u)$, then we decrease by one the labels of the vertices with labels
$\phi(z)+1,\phi(z)+2,\dots,\phi(u)$, and give $z$ new label $\phi(u)$. It is not hard to show that this operation does
not increase the bandwidth of the labeling. We keep repeating this (to each new labeling) until we find no $z$ such
that $\phi(z)<\phi(u)$.

Our heuristic then consists of several (typically one thousand) independent runs, each consisting of three steps:
first randomly ordering the vertices, secondly performing the reverse
Cuthill-McKee algorithm, and thirdly the above improvement procedure. We start each run by a random ordering of
the vertices because the Cuthill-McKee algorithm strongly depends on the initial ordering of vertices
(certainly in the case of symmetric graphs). Our method is thus quite elementary and fast, and as we shall
see in the next section, it gives good results. The interested reader can download the described heuristic algorithm from
\url{https://stuwww.uvt.nl/~sotirovr/research.html}. We note finally that in the literature we did not find
any heuristics for the bandwidth problem that were specifically targeted at graphs with symmetry.

\section{Numerical results} \label{sec:NumR}
In this section we present the numerical results for the bandwidth problem for several graphs. All relaxations were solved
with  SeDuMi \cite{sedumi} using the Yalmip interface \cite{YALMIP} on an Intel Xeon X5680, $3.33$ GHz dual-core
processor with 32 GB memory. To compute orbitals, we used GAP \cite{gap}.

\subsection{Upper bounds}\label{Sect:UpBnd}

The results of our heuristic that we described in Section \ref{Sect:CMheuristic} are given in (some of) the tables in the
next section in the columns named `u.b.'. The computation time for instances with less than 250 nodes (with 1000 runs) is about 30 s.

The obtained output seems to be of good quality: in many cases for which we know the optimal value, this value
is attained; for example for the Hamming graph $H(3,6)$, where we have an upper bound 101, which was shown to be
optimal by Balogh et al.~\cite{baHa}.
We remark that the best obtained upper bound for $H(3,6)$ computed by the Cuthill-McKee algorithm after 1000 random starts,
but without our improvement steps, was 130, while the above described improvement reports the optimal value (101) 22 times (out of the 1000 runs).
Also for the Johnson graphs $J(v,2)$ with $v\in \{4,\ldots,15\}$ the improved reverse Cuthill-McKee heuristic provides sharp bounds.
For more detailed results on the upper bounds of different graphs, see the following section.

\subsection{Lower bounds}

In this section we present several lower bounds on the bandwidth of a graph.
Each such lower bound, i.e., bw$_{\rm eig}$, bw$_{\rm COP}$,  bw$_{\rm QAP}$, and bw$_{\rm fix}$
 is obtained from a lower bound of the corresponding relaxation (OPT$_{\rm eig}$, MC$_{\rm COP}$, MC$_{\rm QAP}$, and MC$_{\rm fix}$, respectively)
 of {\em some} min-cut problem. Indeed, for each graph and each relaxation
we consider several min-cut problems, each corresponding to a different $m$.
In particular, we first computed OPT$_{\rm eig}$ and the corresponding bound for the bandwidth
for all choices of $m=(m_1,m_2,m_3)$ with $m_1\leq m_2$; this can be done in a few seconds for each graph.
The reported bw$_{\rm eig}$ is the best bound obtained in this way, and $m_{\rm eig}$ is the corresponding $m$.
We then computed MC$_{\rm QAP}$ (or  MC$_{\rm COP}$ for the hypercube graph) for all $m$ with $m_1\leq m_2$ and $m_3\geq \overline{m}_3$ (where $\overline{m}=m_{\rm eig}$).
Similarly we computed the bound on the bandwidth that is obtained from  MC$_{\rm fix}$. \\

\subsubsection{The hypercube graph}
The first numerical results, which we present in Table \ref{tab2}, concern the Hamming graph $H(d,2)$,
also known as the hypercube $Q_d$, see Example \ref{ex:hamming}. Because the bandwidth of the hypercube graph was determined
already by Harper \cite{Harper66} as
\[
\sigma_{\infty}(Q_d) = \sum\limits_{i=0}^{d-1} \binom{i}{ \lfloor \frac{i}{2} \rfloor},
\]
this provides a good first test of the quality of our new bound. In the other families of graphs that we tested,
we did not know the bandwidth beforehand, so there the numerical results are really new (as far as we know). Besides that, the tested
hypercube graphs are relatively small, so that we could also compute MC$_{\rm COP}$ (which we do not do in the other examples).

Table \ref{tab2} reads as follows.
The second column contains the number of vertices $n=2^d$ of the graph, and the last column contains the exact values for the bandwidth $\sigma_{\infty}(Q_d)$.
In the third column we give the lower bound on the bandwidth of the graph corresponding to OPT$_{\rm eig}$,
while in the fourth-sixth column the lower bounds correspond  to the solutions of the optimization problems MC$_{\rm COP}$, MC$_{\rm QAP}$, and MC$_{\rm fix}$,
 respectively (obtained as described above).
We remark that when the lower bound bw$_{\rm QAP}$ was tight (this happened for $d=2,3$), we did not compute the bound bw$_{\rm fix}$, since then this is also tight.
Note that for $d=2,3,4$, the latter bound is indeed tight.

\begin{table}[h!]
\begin{center}
\begin{tabular}{|c|c||c|c|c|c||c|}
\hline
$d$ & $n$  & bw$_{\rm eig}$ & bw$_{\rm COP}$  &  bw$_{\rm QAP}$ & bw$_{\rm fix}$ & $\sigma_{\infty}(Q_d)$ \\[1ex] \hline
2 & 4   & 2 & 2 & 2 & 2 &  2 \\[1ex]
3 & 8  & 3 & 4 & 4 &  4 &  4 \\[1ex]
4 & 16  & 4 & 6 & 6 & 7 &  7 \\[1ex]
5 & 32 & 7 & 10 & 10 & 11  &  13 \\\hline
\end{tabular}
\caption{Bounds on the bandwidth of hypercubes $Q_d$.}  \label{tab2}
\end{center}
\end{table}

In  Table \ref{tab2a}, we list the computational times required for solving the optimization problems MC$_{\rm COP}$, MC$_{\rm QAP}$, and MC$_{\rm fix}$
that provide the best bound for the bandwidth.
The computational time to obtain  MC$_{\rm fix}$ for $Q_d$ is equal to the sum of the computational times required for
solving each of the subproblems MC$_{\rm fix}^h$, $h=1,\ldots,d$.
Table \ref{tab2a} provides also the best choice of the vector $m$, i.e., the one that provides the best bandwidth bound for the given optimization problem.
We remark that for $d=3,4$ there are also other such best options for $m$.
Table \ref{tab2b} (see Appendix \ref{ApA}) gives the number of orbitals in the stabilizer subgroups ${\mathcal H}(Q_d(\alpha))$,
for different $\alpha= \{ 1,\ldots, n \} \backslash \{ r_1,r_2 \}$.

\begin{table}[h!]
\begin{center}
\begin{tabular}{|c||c|c|c|c|c|c|c|}
\hline
$d$ & $m_{\rm eig}$ &  MC$_{\rm COP}$ & $m_{\rm COP}$ & MC$_{\rm QAP}$ & $m_{\rm QAP}$ & MC$_{\rm fix}$ & $m_{\rm fix}$ \\[1ex] \hline
2 & $[1,2,1]$  & 0.51 & $[1,2,1]$ & 0.05  & $[1,2,1]$ &  -- & -- \\[1ex]
3 & $[3,3,2]$  & 1.29 & $[2,3,3]$ & 0.53  & $[2,3,3]$ & --  & -- \\[1ex]
4 & $[6, 7,3]$ &23.55  & $[4,7,5]$ &   1.32  & $[4,7,5]$  & 7.12 & $[4,6,6]$   \\[1ex]
5 & $[11,15,6]$ &1019.83  & $[10,14,8]$ & 1.86 & $[10,14,8]$  & 74.56 & $[10,12,10]$  \\\hline
\end{tabular}
\caption{$Q_d$: time (s) to solve relaxations  and corresponding $m$.}  \label{tab2a}
\end{center}
\end{table}

\subsubsection{The Hamming graph}
Next, we give lower and upper bounds on the bandwidth of the Hamming graph $H(d,q)$ for $q>2$, see  Example \ref{ex:hamming}.
Because the two-dimensional Hamming graph $H(2,q)$ (also known as the lattice graph) has bandwidth equal to $\frac{(q+1)q}{2}-1$ (see \cite{HenSti92}),
we computed bounds on the bandwidth for the next two interesting groups of Hamming graphs, i.e., $H(3,q)$ and $H(4,q)$.
Tables \ref{tab4}, \ref{tab4a}, and \ref{tab4b} are set-up similarly as in the case of the hypercube graphs, except for the last column of
Table \ref{tab4}, where instead of the exact value of the bandwidth we now give the upper bound obtained by the heuristic as described in Section \ref{Sect:UpBnd}.

Concerning this upper bound, we remark that it matches the exact bandwidth for $H(3,6)$ (see \cite{baHa}). Besides this case,
the upper bounds for the Hamming graph from the literature are weaker than ours.
For example, Harper's  upper bound \cite[Claim 1]{Harper03}
for $H(3,4)$ is $33$, while the upper bound for $H(3,5)$ by Berger-Wolf and Reingold \cite[Thm. 2]{BergWoRe:02} is $84$.

For the lower bounds, it is interesting to note that bw$_{\rm eig}$ equals bw$_{\rm QAP}$ for $H(3,4)$,
while the lower bound bw$_{\rm fix}$ is (strictly) the best for all graphs in the table.

Concerning Table \ref{tab4a}, we note that for all instances there are more options for $m$ that provide the best bound.
The number of orbitals in ${\mathcal H}([H(4,3)](\alpha))$ for each subproblem as given in Table \ref{tab4b} is large,
but we are able to compute MC$_{\rm fix}$ because the adjacency matrix $A(\alpha)$ has order only $79$.

\begin{table}[h!]
\begin{center}
\begin{tabular}{|c|c|c||c|c|c||c|}
\hline
 $d$ & $q$ & $n$ &  bw$_{\rm eig}$ &  bw$_{\rm QAP}$ & bw$_{\rm fix}$ & u.b. \\[1ex] \hline
3&3 &  27   & 9 & 10  &  12  & 13 \\[1ex]
3& 4 & 64    & 22 & 22 & 25 & 31  \\[1ex]
3&5 & 125   & 42 & 43 & 47 & 60 \\[1ex]
3&6 & 216  & 72 & 74  & 78 &  101\\\hline
%4&2 &  16 & 5 & 6 & 7 & 7 \\
4&3 &  81   & 21 & 23 &  26  &  35 \\ \hline
\end{tabular}
\caption{Bounds on the bandwidth of the Hamming graphs $H(3,q)$ and $H(4,q)$.}  \label{tab4}
\end{center}
\end{table}

\begin{table}[h!]
\begin{center}
\begin{tabular}{|c|c||c|c|c|c|c|}
\hline
$d$& $q$ & $m_{eig}$ & MC$_{\rm QAP}$ &  $m_{\rm QAP}$ & MC$_{\rm fix}$  & $m_{\rm fix}$\\[1ex] \hline
3&3 &  $[9,10,8]$ &  0.29  & $[9,10,8]$ &  44.03 & $[4,12,11]$  \\[1ex]
3&4 &  $[21,22,21]$ &  2.80 & $[21,22,21]$ & 176.38 & $[20,20,24]$ \\[1ex]
3&5 &  $[35,41,49]$& 15.50 & $[42,43,40]$ &  536.65 & $[38,41,46]$\\[1ex]
3&6 &  $[61,71,84]$ & 76.20 &$[72,74,70]$ & 1756.88 & $[63,76,77]$ \\\hline
%4&2 &  $[6,7,3]$ & 1.32 & $[6,7,3]$ &  7.96 &  $[3,7,6]$ \\[1ex]
4&3 &  $[20,27,34]$ &  9.10  & $[29,30,22]$ & 5877.33  & $[22,34,25]$  \\\hline
\end{tabular}
\caption{ $H(3,q)$ and $H(4,q)$:  times (s) to solve relaxations and corresponding $m$.}  \label{tab4a}
\end{center}
\end{table}

\subsubsection{The $3$-dimensional generalized Hamming graph}
In Tables \ref{tab5}, \ref{tab5a}, and \ref{tab5b}, we present the analogous numerical results for the
$3$-dimensional generalized Hamming graph $H_{q_1,q_2,q_3}$, as defined in Example \ref{ex:3dimgenham}.
To the best of our knowledge there are no other lower and/or upper bounds for the bandwidth of $H_{q_1,q_2,q_3}$ in the literature.
Our results show that the new bound can be significantly better that the eigenvalue bound.
For instance, the eigenvalue lower bound on the bandwidth of $H_{3,4,5}$ is 16, while bw$_{\rm fix}$ is 24.\\

\begin{table}[h!]
\begin{center}
\begin{tabular}{|c|c|c|c||c|c|c||c|}
\hline
$q_1$ & $q_2$ & $q_3$ & $n$ &  bw$_{\rm eig}$ &  bw$_{\rm QAP}$ & bw$_{\rm fix}$ & u.b. \\[1ex] \hline
 2& 3 & 3 & 18   &  5 & 8  & 9   & 9 \\[1ex]
 2 & 3 & 4 & 24 & 6 & 10  &  11  & 12 \\[1ex]
 2 & 3 & 5 & 30 & 6 & 11  &  13  & 15 \\[1ex]
 2 & 4 & 4 & 32 & 7 & 12 &  14 & 16 \\[1ex]
 3 & 3 & 4 & 36 & 11 & 13 &  15 & 17 \\[1ex]
 3 & 3 & 5 & 45 & 13 & 16 & 19  & 21 \\[1ex]
 3 & 4 & 4 & 48 & 14 & 17 & 20  & 23 \\[1ex]
 3 & 4 & 5 & 60 & 15 & 21 & 24  & 29 \\[1ex]
 \hline
\end{tabular}
\caption{Bounds on the bandwidth of $H_{q_1,q_2,q_3}$.}  \label{tab5}
\end{center}
\end{table}

\begin{table}[h!]
\begin{center}
\begin{tabular}{|c|c|c||c|c|c|c|c|}
\hline
$q_1$& $q_2$ & $q_3$ & $m_{eig}$ & MC$_{\rm QAP}$ &  $m_{\rm QAP}$ & MC$_{\rm fix}$  & $m_{\rm fix}$\\[1ex] \hline
2& 3 & 3    & $[6,8,4]$ & 1.27 & $[4,8,6]$ & 45.94 & $[5,5,8]$ \\[1ex]
2 & 3 & 4 & $[8,11,5]$ & 1.11 & $[6,10,8]$ & 209.49 & $[6,8,10]$\\[1ex]
2& 3 & 5 & $[11,14,5]$ & 2.20 & $[8,13,9]$ & 240.93 & $[8,10,12]$ \\[1ex]
2& 4 & 4 & $[12,14,6]$ & 2.02 & $[10, 12,10]$ & 99.68 & $[9,10,13]$ \\[1ex]
3& 3 & 4 & $[12,14,10]$ & 1.19 & $[6,18,12]$ & 558.94 & $[5,17,14]$ \\[1ex]
3& 3 & 5 & $[14, 19, 12]$ & 4.35 & $[15,16,14]$ & 525.82 & $[13,14,18]$ \\[1ex]
3& 4 & 4 & $[16,19,13]$ & 3.56 & $[14,19,15]$ & 693.19 & $[14,15,19]$\\[1ex]
3& 4 & 5 & $[18,28,14]$ & 6.37 & $[20,21,19]$ & 3702.17& $[18,19,23]$\\[1ex]
 \hline
\end{tabular}
\caption{$H_{q_1,q_2,q_3}$: time (s) to solve relaxations and corresponding $m$.}  \label{tab5a}
\end{center}
\end{table}

\subsubsection{The Johnson graph}
In Tables \ref{tab3}, \ref{tab3a}, and \ref{tab3b}, we present the analogous numerical results for the Johnson graph $J(v,d)$, as defined in Example \ref{ex:johnson}.
In particular, we provide bounds for $J(v,3)$, with $v\in \{6,\ldots, 11\}$, and $J(8,4)$.
The bandwidth of the Johnson graph $J(v,2)$ (also known as the triangular graph) has been determined by Hwang and Lagarias \cite{HwLa77},
and equals $\lfloor v^2/4 \rfloor + \lceil v/2 \rceil -2$.
We remark that when the bound bw$_{\rm QAP}$ was tight (this happened for $v=6,7$), we did not compute the bound bw$_{\rm fix}$, since then the latter is also tight.
Indeed, it thus follows that the bandwidth of $J(6,3)$ equals $13$ and the bandwidth of  $J(7,3)$ equals $22$.

We also remark here that there are more options for $m$ that provide the best bound, in particular for $J(v,3)$, with $v\geq 8$.
For example, the lower bound for the bandwidth of $J(8,3)$ is equal to $31$ for the vectors $[12,14,30]$ and $[13,13,30]$.
Table \ref{tab3b}, see Appendix \ref{ApA}, provides the number of orbitals from the stabilizer subgroups
${\mathcal H}([J(v,d)](\alpha))$,  $\alpha= \{ 1,\ldots, n \} \backslash \{ r_1,r_2 \}$ for $d=3$ and $d=4$, respectively.
Since for $d=4$ the number of orbitals increases significantly when $v$ increases from eight  to nine,
we could not compute the new lower bound bw$_{\rm fix}(J(9,4))$.
However, we obtained that bw$_{\rm eig}(J(9,4))=49$, and bw$_{\rm QAP}(J(9,4))=52$ in 23.12 seconds. \\

\begin{table}[h!]
\begin{center}
\begin{tabular}{|c|c|c||c|c|c||c|}
\hline
$v$ & $d$ & $n$ & bw$_{\rm eig}$ &  bw$_{\rm QAP}$ & bw$_{\rm fix}$ & u.b. \\[1ex] \hline
6 & 3& 20  & 10 & 13 &  13 & 13  \\[1ex]
7 & 3& 35  & 17 & 22 & 22 & 22 \\[1ex]
8 & 3& 56  & 25 & 29 & 31 & 34 \\[1ex]
9 & 3& 84  & 36 &  40 & 43 & 49 \\[1ex]
10 & 3&120  & 50 & 53 & 57 & 68\\[1ex]
11 & 3&165  & 68 & 69 & 74 & 92\\[1ex] \hline
8 & 4& 70  & 28 & 33 & 37 & 40 \\ \hline % \\[1ex]
\end{tabular}
\caption{Bounds on the bandwidth of $J(v,3)$  and $J(v,4)$.}  \label{tab3}
\end{center}
\end{table}

\begin{table}[h!]
\begin{center}
\begin{tabular}{|c|c||c|c|c|c|c|}
\hline
$v$ & $d$& $m_{\rm eig}$ & MC$_{\rm QAP}$ &  $m_{\rm QAP}$ & MC$_{\rm fix}$  & $m_{\rm fix}$\\[1ex] \hline
6 &3& $[5,6,9] $ & 0.26  & $[3,5,12]$ & --  & --   \\[1ex]
7 &3& $[8,11,16]$ &  0.87  & $[7,8,20]$ & -- %57.87
& -- \\[1ex]
8 &3& $[16, 16, 24]$ & 2.26  & $[11,18, 27]$ & 194.24  & $[13,13,30]$ \\[1ex]
9 & 3& $[22,27,35]$ & 5.57  & $[14,32,38]$ & 558.01 & $[15,27,42]$ \\[1ex]
10 & 3& $[ 28,43,49]$ & 14.72 & $[18,51,51]$ & 865.89 & $[32,32,56]$ \\[1ex]
11 & 3& $[45,53,67]$ &  34.94 & $[43,57,65]$ & 1607.29  &  $[38,56,71]$ \\[1ex] \hline
8 &4 & $[16,27,27]$ & 2.31 & $[13, 26,31]$  &  368.49 & $[12,22,36]$   \\ \hline %\\[1ex]
\end{tabular}
\caption{$J(v,3)$ and $J(v,4)$: time (s) to solve relaxations and corresponding $m$.}  \label{tab3a}
\end{center}
\end{table}

\subsubsection{The Kneser graph}
In Tables \ref{tab7} and \ref{tab7a}, we finally present our bounds and computational times for the Kneser
graph $K(v,2)$  with $v\in \{5,\ldots, 8\}$ and $K(v,3)$  with $v\in \{7,\ldots, 10\}$, see Example \ref{ex:johnson}.
We remark that the orbitals of aut($K(v,d))$ and aut($J(v,d)$) are the same, and the corresponding orbitals
in the stabilizer subgroups for these two graphs are also the same.

We note that Juvan and Mohar \cite{JuMo:99} presented  general lower and upper bounds for the Kneser graph.
Their lower bound however is an eigenvalue bound that is weaker than the eigenvalue bound bw$_{\rm eig}$ by Helmberg et al.~\cite{HelReMoPo:95}.
Also their general upper bound is (much) weaker than our computational results.

Since the bound bw$_{\rm QAP}$ for the Petersen graph $K(5,2)$ is tight, we did not compute  bw$_{\rm fix}$ for this graph.
Note that it is a folklore result that the bandwidth of the Petersen graph equals $5$.

\begin{table}[h!]
\begin{center}
\begin{tabular}{|c|c|c||c|c|c||c|}
\hline
$v$ & $d$ & $n$ & bw$_{\rm eig}$ &  bw$_{\rm QAP}$ & bw$_{\rm fix}$ & u.b. \\[1ex] \hline
5 & 2 & 10  & 4 & 5 & 5 & 5 \\[1ex]
6 & 2 & 15  & 9 & 9 & 10  & 10 \\[1ex]
7 & 2 & 21  & 14 & 14 & 15  & 16 \\[1ex]
8 & 2 & 28  & 20 & 20 &  22 & 23 \\[1ex]
\hline
7 & 3 & 35  & 10 & 12 & 12  & 15 \\[1ex]
8 & 3 & 56  & 25 & 26 & 27  & 33 \\[1ex]
9 & 3 & 84  & 45 & 47 & 48  & 59 \\[1ex]
10 & 3 & 120  & 72 & 75 & 76 & 90 \\[1ex]
\hline
\end{tabular}
\caption{Bounds on the bandwidth of $K(v,2)$ and $K(v,3)$.}  \label{tab7}
\end{center}
\end{table}

\begin{table}[h!]
\begin{center}
\begin{tabular}{|c|c||c|c|c|c|c|}
\hline
$v$ & $d$& $m_{\rm eig}$ & MC$_{\rm QAP}$ &  $m_{\rm QAP}$ & MC$_{\rm fix}$  & $m_{\rm fix}$\\[1ex] \hline
5 &2 & $[3,4,3]$ & 0.45  & $[3,4,3]$ & --  & --   \\[1ex]
6 &2 & $[3,4,8]$ &  0.43   & $[3,4,8]$ & 2.18 & $[3,3,9]$   \\[1ex]
7 &2 & $[4,4,13]$ &  0.49  & $[4,4,13]$ & 3.71 & $[3,4,14]$   \\[1ex]
8 &2 & $[3,6,19]$ &  0.80  & $[3,6,19]$ & 6.75 & $[4,4,20]$ \\[1ex]
\hline
7 & 3 & $[12, 14, 9]$ &  1.77  & $[11, 14, 10]$ & 60.81 & $[11, 14, 10]$   \\[1ex]
8 & 3 & $[13,19,24]$ &  2.29  & $[15,18,23]$ & 180.09 & $[14,16,26]$   \\[1ex]
9 & 3 & $[16,24,44]$ &  6.74  & $[19,22,43]$ & 561.12 & $[21,16,47]$   \\[1ex]
10 & 3 & $[19,30,71]$ &   14.09 & $[24,26,70]$ & 1043.87 &  $[25,21,74]$  \\[1ex]
\hline
\end{tabular}
\caption{$K(v,2)$ and  $K(v,3)$: time (s) to solve relaxations  and corresponding $m$.}  \label{tab7a}
\end{center}
\end{table}

\appendix

\section{$*$-isomorphism from Example \ref{ex:Kkl}} \label{Appendix*isomorph}

The associated $*$-isomorphism $\varphi$ satisfies:
\[
\tiny \varphi(B_1) =
\begin{pmatrix}
1 & & & \\
  & 0 & & \\
  &  & 0 & \\
  && & 1 & 0 &0 \\
 & && 0 & 0 & 0 \\
 & && 0 & 0 & 0 \\
 \end{pmatrix}, \; \quad
\varphi(B_2) =
\begin{pmatrix}
-1 & & & \\
  & 0 & & \\
  &  & 0 & \\
  && & m_1-1 & 0 &0 \\
 & && 0 & 0 & 0 \\
 & && 0 & 0 & 0 \\
 \end{pmatrix},
\]
\[
\tiny
 \varphi(B_3) =
 \sqrt{m_1 m_2}
\begin{pmatrix}
0 & & & \\
  & 0 & & \\
  &  & 0 & \\
  && & 0 & 1 &0 \\
 & && 0 & 0 & 0 \\
 & && 0 & 0 & 0 \\
 \end{pmatrix}, \;\quad
\varphi(B_4) =
 \sqrt{m_1m_3}
\begin{pmatrix}
0 & & & \\
  & 0 & & \\
  &  & 0 & \\
  && & 0 & 0 & 1  \\
 & && 0 & 0 & 0 \\
 & && 0 & 0 & 0 \\
 \end{pmatrix}, \; \quad
\]
\[
\tiny
 \varphi(B_5) =
 \sqrt{m_1 m_2}
\begin{pmatrix}
0 & & & \\
  & 0 & & \\
  &  & 0 & \\
  && & 0 & 0 &0 \\
 & && 1 & 0 & 0 \\
 & && 0 & 0 & 0 \\
 \end{pmatrix}, \;\quad
\varphi(B_6) =
\begin{pmatrix}
0 & & & \\
  & 1 & & \\
  &  & 0 & \\
  && & 0 & 0 & 0  \\
 & && 0 & 1 & 0 \\
 & && 0 & 0 & 0 \\
 \end{pmatrix}, \;\quad
\]
\[
\tiny
 \varphi(B_7) =
\begin{pmatrix}
0 & & & \\
  & -1 & & \\
  &  & 0 & \\
  && & 0 & 0 &0 \\
 & && 0 & m_2-1 & 0 \\
 & && 0 & 0 & 0 \\
 \end{pmatrix}, \;\quad
\varphi(B_8) =
 \sqrt{m_2 m_3}
\begin{pmatrix}
0 & & & \\
  & 0 & & \\
  &  & 0 & \\
  && & 0 & 0 & 0  \\
 & && 0 & 0 & 1 \\
 & && 0 & 0 & 0 \\
 \end{pmatrix}, \;\quad
\]
\[
\tiny
 \varphi(B_9) =
 \sqrt{m_1 m_3}
\begin{pmatrix}
0 & & & \\
  & 0 & & \\
  &  & 0 & \\
  && & 0 & 0 &0 \\
 & && 0 & 0 & 0 \\
 & && 1 & 0 & 0 \\
 \end{pmatrix}, \;\quad
\varphi(B_{10}) =
 \sqrt{m_2 m_3}
\begin{pmatrix}
0 & & & \\
  & 0 & & \\
  &  & 0 & \\
  && & 0 & 0 & 0  \\
 & && 0 & 0 & 0 \\
 & && 0 & 1 & 0 \\
 \end{pmatrix}, \;\quad
\]
\[
\tiny
 \varphi(B_{11}) =
 \begin{pmatrix}
0 & & & \\
  & 0 & & \\
  &  & 1 & \\
  && & 0 & 0 &0 \\
 & && 0 & 0 & 0 \\
 & && 0 & 0 & 1 \\
 \end{pmatrix}, \;\quad
\varphi(B_{12}) =
\begin{pmatrix}
-1 & & & \\
  & 0 & & \\
  &  & 0 & \\
  && & 0 & 0 & 0  \\
 & && 0 & 0 & 0 \\
 & && 0 & 0 & m_3-1 \\
 \end{pmatrix}.
\]

\section{Symmetry reduction of $\mbox{MC}_{\rm fix}^h$} \label{AppendixMcfixh}

Let $G$ be an undirected graph on $n$ vertices with adjacency matrix $A$ and $t$ orbitals $\mathcal{O}_h$ ($h=1,2,\ldots,t$).
Let $(s_1,s_2)$ be an arbitrary edge in the  cut graph $G_{m_1,m_2,m_3}$ with the adjacency matrix $B$,
and $(r_{h1},r_{h2})$ be an arbitrary pair of vertices in $\mathcal{O}_h$ ($h=1,2,\ldots,t$).
We let $\alpha^h = \{ 1,\ldots, n \}  \backslash \{ r_{h1}, r_{h2} \}$ and $\beta =\{ 1,\ldots, n \}  \backslash \{ s_1,s_2 \}$.
Now, the  relaxation $\mbox{MC}_{\rm fix}^h$ (see page \pageref{Mcfixhh}) reduces to
\begin{equation}\label{reducedEQP}
\begin{array}{rl}
\min  & \frac{1}{2} \sum\limits_{i=1}^d ~\sum\limits_{j=1}^{12} p_i^{-1}  \tr(A(\alpha^h)A_i) x^{(i)}_j  \\[3ex]
& + \sum\limits_{i\in {\mathcal I}_{\mathcal A}} \sum\limits_{j\in \{1,6,11\}} (q_j p_i)^{-1} B(\beta,s_1)^{\mathrm{T}}\diag(B_j) A(\alpha^h,r_{h1})^{\mathrm{T}}\diag(A_i)x^{(i)}_j     \\[3ex]
 & ~~~+ \sum\limits_{i\in {\mathcal I}_{\mathcal A}} \sum\limits_{j\in \{1,6,11\}} (q_j p_i)^{-1} B(\beta,s_2)^{\mathrm{T}} \diag(B_j)A(\alpha^h,r_{h2})^{\mathrm{T}}\diag(A_i) x^{(i)}_j
 +\frac{1}{2}d^h  \\[3ex]
{\rm s.t.} &  \sum\limits_{i\in {\mathcal I}_{\mathcal A}} x^{(i)}_1=q_1,
 \quad \sum\limits_{i\in {\mathcal I}_{\mathcal A}} x^{(i)}_6=q_6
\quad \sum\limits_{i\in {\mathcal I}_{\mathcal A}} x^{(i)}_{11}=q_{11}   \\[3ex]
%& \sum\limits_{j=1}^d x^{(i)}_j = q_i, \quad i=1,\ldots,12 \\[2ex]
& \sum\limits_{i=1}^d ~\sum\limits_{j=1}^{12} q_j^{-1}x^{(i)}_j B_j = J_{n-2} \\
& \sum\limits_{j=1}^{12} x^{(i)}_j = p_i, \quad i=1,\ldots,d  \\[2ex]
& \sum\limits_{i=1}^{d} \sum\limits_{j=1}^{12}  \frac{1}{q_j p_i} x^{(i)}_j (B_j\otimes A_i) \succeq 0\\[2.5ex]
& x^{(i)}_j \geq 0, ~ x^{(i)}_{j^*} =x^{(i^{*})}_j,
\quad i=1,\ldots,d, ~j=1,\ldots,12,
\end{array}
\end{equation}
where $B_j$ ($j=1,\ldots, 12$) is defined in Example  \ref{ex:Kkl}, and $\{A_i:i=1,\ldots,d\}$ spans  ${\mathcal H}(A(\alpha^h))$.
The set ${\mathcal I}_{\mathcal A} := {\mathcal I}_{{\mathcal H}(A(\alpha^h))}$ is as in Definition \ref{def:coherent config},
$p_i = \tr(J_{n-2}A_i)$, $i=1,\ldots, d$,  $q_j = \tr(J_{n-2}B_j)$, $j=1,\ldots, 12$.
The constraint $x^{(i)}_{j^*} =x^{(i^{*})}_j$ requires that the variables $x^{(i)}_{j}$ form complementary pairs.
The SDP relaxation (\ref{reducedEQP}) can be further simplified by exploiting the $*$-isomorphism associated
to ${\mathcal H}(B(\beta))$, see  Appendix \ref{Appendix*isomorph}.

\section{Orbitals in stabilizer subgroups} \label{ApA}

\begin{table}[h!]
\begin{center}
\begin{tabular}{|c||ccccc|}
\hline
 $d$ &   \multicolumn{5}{c}{$\sharp$ orbitals} \vline\\\cline{1-6}\hline\hline
%% $d$ & (1,2) & (1,4) & (1,8) & (1,16) & (1,32)  \\[1ex] \hline
4  & 80   & 100 & 80 & 35 &  --\\
5  & 140 & 200 & 200 & 140 & 56 \\\hline
\end{tabular}
\caption{Number of orbitals in  ${\mathcal H}(Q_d(\alpha))$. }  \label{tab2b}
\end{center}
\end{table}

\begin{table}[h!]
\begin{center}
\begin{tabular}{|c|c||cccc|}
\hline
$d$ & $q$ &  \multicolumn{4}{c}{$\sharp$ orbitals} \vline\\ %\cline{3-6}
\hline\hline
3 &3 & $135$ & 225 & 165 & --\\
3 & 4  & $150$ & 275 & 220 & -- \\
3& 5 & $150$ & 275 & 220 & -- \\
3 &6 & $150$ & 275 & 220& -- \\\hline
%4 & 2 & $80$ & 100 & 80  & 35 \\
4 &3 & $315$ & 675 & 825  & 495  \\\hline
\end{tabular}
\caption{Number of orbitals  in  ${\mathcal H}([H(d,q)](\alpha))$. }  \label{tab4b}
\end{center}
\end{table}

\begin{table}[h!]
\begin{center}
\begin{tabular}{|c|c|c||ccccccc|}
\hline
$q_1$ & $q_2$ & $q_3$ &  \multicolumn{7}{c}{ $\sharp$ orbitals}  \vline\\[1ex] \hline
 2 & 3 & 3  & 180  & 180 & 60 & 180 & 180 & -- & --  \\
 2 & 3 & 4  & 200  & 180 & 360 & 100 & 200 &180 & 360\\
 2 & 3 & 5 & 200  & 180 & 360 & 100 & 200 & 180 & 360 \\
 2 & 4 & 4 & 200  & 220 &  60 &  200 & 220 & -- & --   \\
 3 & 3 & 4 &  150 & 225  & 450 & 225 & 450 & -- & --  \\
 3 & 3 & 5 & 150  & 225 & 450 & 225 & 450 & -- & --   \\
 3 & 4 & 4 &   250 & 275 & 135 & 450 & 495 &-- & --  \\
 3 & 4 & 5 & 250  & 250 & 500 & 225 & 450 & 450 & 900 \\
 \hline
\end{tabular}
\caption{Number of orbitals  in  ${\mathcal H}(H_{q_1,q_2,q_3}(\alpha))$}.  \label{tab5b}
\end{center}
\end{table}

\begin{table}[h!]
\begin{center}
\begin{tabular}{|c|c||cccc|}
\hline
$v$ & $d$&  \multicolumn{4}{c}{ $\sharp$ orbitals} \vline\\\cline{3-5}\hline\hline
6 & 3 & 88  & 88 & 24 &  --\\
7 & 3 & 195 & 257 &  90 & --  \\
8 & 3 & 220  & 333 & 158 &--  \\
9 & 3 & 227 & 361 & 203  &--\\
10 & 3 & 228 &  368 & 220 &--\\
11 & 3 & 228 & 369 & 225 &-- \\\hline
8 & 4 & 220 & 358 &  220 &  46 \\
9 & 4 &  484 & 916 & 742 & 195 \\ \hline
\end{tabular}
\caption{Number of orbitals in  ${\mathcal H}([J(v,3)](\alpha))$ and ${\mathcal H}([J(v,4)](\alpha))$.}  \label{tab3b}
\end{center}
\end{table}

\end{document}